\documentstyle[12pt]{article}
\input {cyracc.def}
\font\ququ=cmr10 scaled \magstep1
\font\tencyr=wncyr7 scaled \magstep1
\def\rus{\tencyr\cyracc}

\newcommand{\re}[1]{\mbox{\rm  (\ref{#1})}}

\newenvironment{proof}
{\noindent {\sl Proof.}\quad }{\hfill 
$\square$ \vskip1.1ex\noindent }

\newenvironment{proof*}
{\noindent {\sl Proof.}\quad }{\hfill 
$\square$}

\renewcommand{\theequation}{\thesection .\arabic{equation}}

\catcode`\@=\active
\catcode`\@=11
\def\@eqnnum{\hbox to .01pt{}\rlap{\bf \hskip -\displaywidth(\theequation)}}
\catcode`\@=12

\newenvironment{s}[1]
{ \vskip1.2ex \refstepcounter{equation}
\noindent {\bf \theequation\quad #1.} \begin{sl}}{\end{sl}
\vskip1.1ex\noindent }

\newenvironment{rem}[1]
{ \vskip1.2ex \refstepcounter{equation}
\noindent {\bf \theequation\quad {#1}.} }{ \vskip1.1ex\noindent }

\newenvironment{subs}[1]
{\vskip1.2ex \refstepcounter{equation}
\noindent {\bf (\theequation)\quad #1.} }{\quad}

\newcommand {\sekt}[1]
{{\vskip2.5ex\refstepcounter{section}\setcounter{equation}{0}
\noindent\large \bf \thesection\quad #1
\vskip1.5ex\noindent}}

\newcommand {\ah}{{\frak a}}

\newcommand {\ce}{{\frak c}}
\newcommand {\fe}{{\frak E}}

\newcommand {\g}{{\frak g}}

\newcommand {\ka}{{\frak k}}
\newcommand {\el}{{\frak l}}
\newcommand {\n}{{\frak n}}

\newcommand {\p}{{\frak p}}

\newcommand {\te}{{\frak t}}
\newcommand {\tri}{{\frak sl}_2}

\newcommand {\z}{{\frak z}}
\newcommand {\esi}{\varepsilon}
\newcommand {\ap}{\alpha}

\newcommand {\lb}{\lambda}

\newcommand {\tih}{\tilde h}

\newcommand {\V}{{\Bbb V}}

\newcommand {\ca}{{\cal A}}
\newcommand {\cb}{{\cal B}}

\newcommand {\N}{{\cal N}}
\newcommand {\co}{{\cal O}}
\newcommand {\cs}{{\cal S}}

\newcommand {\cw}{{\cal W}}
\newcommand {\hb}{{\bf h}}

\newcommand {\eb}{{\bf{e}}}
\newcommand {\db}{{\bf{d}}}

\newcommand {\ad}{{\mathrm{ad\,}}}
\newcommand {\Ad}{{\mathrm{Ad\,}}}

\newcommand {\rk}{{\mathrm{rk\,}}}

\newcommand {\GR}[2]{{\mathrm{{\bf #1}}}_{#2}}
\newcommand {\GRt}[2]{{\tilde{\mathrm{{\bf #1}}}}_{#2}}
\newcommand {\ov}{\overline}
\newcommand {\un}{\underline}

\newcommand {\vno}[1]{\vskip#1 ex\noindent}
\newcommand {\rar}{\rightarrow}

\newcommand {\qu}{\hfill $\square$}
\newcommand {\beq}{\begin{equation}}
\newcommand {\eeq}{\end{equation}}
\newcommand {\pn}{{\it pn}-}
\newcommand {\spr}{{\it spr\/}-}
\newcommand {\rn}{regular nilpotent }

\newfam\Bbbfam\newfam\eufam%
\font\Bbbfont=msbm10 scaled 1200%
\font\olala=msam10 scaled 1200%
\font\frak=eufm10 scaled 1400%
\font\Bbbsmallfont=msbm10%
\textfont\Bbbfam=\Bbbfont\scriptfont\Bbbfam=\Bbbsmallfont
\font\euzw=eufm10 scaled 1200%
\font\euac=eufm7 scaled 1200%
\textfont\eufam=\euzw\scriptfont\eufam=\euac%
\def\frak{\fam\eufam}%
\def\Bbb{\fam\Bbbfam}%

\def\varnothing{\hbox {\Bbbfont\char'077}}
\def\square{\hbox {\olala\char"03}}
\def\Bbbk{\hbox {\Bbbfont\char'174}}
\def\cyeq{\hbox {\olala\char'064}}

\begin{document}
\setlength{\parskip}{3pt plus 5pt minus 0pt}
\vskip1ex

\noindent
{\Large \bf Nilpotent pairs, dual pairs, and sheets}
\bigskip \\
{\bf Dmitri I. Panyushev}

\medskip
\smallskip

\noindent{\large \bf Introduction}
\vno{2}%
Recently, {\sc V.\,Ginzburg} introduced the
notion of a {\it principal nilpotent pair} (=~\pn pair) 
in a semisimple Lie algebra $\g$ \cite{vitya}.
It is a double counterpart of the notion of a regular nilpotent
element in $\g$. A pair $\eb=(e_1,e_2)\in\g\times\g$ is called 
{\it nilpotent}, if
$[e_1,e_2]=0$ and there exists a pair $\hb=(h_1,h_2)$ of semisimple elements
such that $[h_1,h_2]=0$, $[h_i,e_j]=\delta_{ij}e_j$ $(i,j\in\{1,2\})$.
A \pn pair $\eb$ is a nilpotent pair
such that the simultaneous
centralizer $\z_\g(\eb)$ is of dimension $\rk\g$. Note that, by a theorem of
{\sc R.\,Richardson},  $\rk\g$ is the least possible value for this 
dimension. Evident similarity
between the ``double"  and ``ordinary" theory is manifestly seen in the 
following results of \cite{vitya}: $\z_\g(\hb)$ is a Cartan subalgebra;
the eigenvalues of $\ad h_1$ and $\ad h_2$ are integral; both $e_1$ and $e_2$
are Richardson elements; $Z_G(\eb)$ is a
connected Abelian unipotent group; $Z_G(\eb)$ acts transitively  on the set
of semisimple pairs satisfying the above commutator relations.
Excerpts from Ginzburg's theory, which by no means exhaust \cite{vitya},
are presented in section~1.
\par 
The aim of this article is to develop the theory of nilpotent pairs a bit
further and to present some applications of it to dual pairs and sheets in
semisimple Lie algebras. In section 2, 
it is shown that a considerable part of the above-mentioned results
can be extended to the nilpotent pairs with
$\dim\z_\g(e_1,e_2)=\rk\g+1$. Such pairs are called {\it almost\/} \pn pairs.
Although almost \pn pairs share many properties with \pn pairs, with
similar proofs, some new phenomena do occur for the former. For instance,
it is shown that the totality of almost \pn pairs breaks into two natural 
classes \re{alm4}.
One of the distinctions between them is that the eigenvalues of $\ad h_i$
($i=1,2$) are integral for the first class and non-integral for the second 
class. We also give a description of $Z_G(\eb)$ for both classes.
It is worth noting that the very existence of almost \pn pairs is
a purely "double" phenomenon, because the dimension of "ordinary" orbits is 
always even.
\par 
It is not always the case that a nilpotent pair can be embedded 
in $\tri\oplus\tri$. The pairs admitting such an embedding are called 
{\it rectangular}. Then, as usual, powerful $\tri$-machinery invented by
{\sc V.V.\,Morozov} and {\sc E.B.\,Dynkin} in 40's makes the life much easier.
For instance, a complete classification of
rectangular \pn pairs is found in \cite{wir}, while a classification
of arbitrary \pn pairs is not yet known. Some results on rectangular pairs,
in particular almost principal ones, are presented in section~3.
\\[.5ex]
Section 4 establishes a relationship between nilpotent pairs and dual pairs.
Given a quadruple $(\eb,\hb)$ satisfying the 
commutator relations as above, we show that
$\ka_1=\z_\g(e_1,h_1)$ and $\ka_2=\z_\g(e_2,h_2)$ form a dual pair in $\g$
under certain constraints (see \ref{dual-main}).
These constraints are satisfied by the \pn pairs and almost \pn pairs.
In the principal or almost principal case, this dual pair is
reductive if and only if $\eb$ is rectangular. Moreover, if $\eb$ is a 
rectangular \pn pair, then $(\ka_1,\ka_2)$ is $S$-irreducible in the sense of
{\sc H.\,Rubenthaler}.
\par 
In section 5, we describe another class of rectangular pairs such that
$(\ka_1,\ka_2)$ is a dual pair. These pairs are called
{\it semi-principal\/} rectangular. 
It is worth to note that, since $\ka_1$ is a centralizer, 
$(\ka_1,\z_\g(\ka_1))$ is a dual pair. So, the point is that 
$\ka_2=\z_\g(\ka_1)$
comes up also as centralizer attached to the second member of the pair.
\par 
As a by-product of our study of semi-principal pairs, we found that 
the double centralizer of some $\tri$-triples  has beautiful properties.
It turns out that this phenomenon, appropriately formalized, has 
some application to sheets. Let $\{e,\tih,f\}$ be an $\tri$-triple.
Both the triple and $e$ are called {\it excellent}, if $e$ is even
and $\dim\z_\g(\z_\g(\tih))=
\rk\z_\g(\z_\g(e,\tih,f))$.
In section 6, we show that the excellent triples enjoy the following
properties: $\z_\g(\z_\g(e,\tih,f))$ is semisimple; $\z_\g(\z_\g(e))$ (resp.
$\z_\g(\z_\g(\tih))$)
is the centralizer of $e$ (resp. $\tih$) in $\z_\g(\z_\g(e,\tih,f))$. Then we 
consider the sheet $\cs$ associated to $\{e,\tih,f\}$.
It is proven that $\cs$ is smooth and has a section,
which is an affine space, and that it is the {\it unique} sheet
containing $e$, see \re{section}.
This applies, in particular, to both members of rectangular \pn pairs.
\par 
In section 7, we classify the excellent elements in the simple Lie
algebras. 
\vno{.5}%
The ground field $\Bbbk$ is algebraically closed and of 
characteristic zero. Throughout, 
$\g$ is a semisimple Lie algebra and $G$ is its adjoint group.
For any set $M\subset\g$, let $\z_\g(M)$ (resp. $Z_G(M)$) denote 
the centralizer of $M$ in $\g$ (resp. in $G$). 
For $M=\{a,\dots,z\}$, we simply write $\z_g(a,\dots,z)$ or $Z_G(a,\dots,z)$.
If $N\subset G$, then $Z_G(N)$ stands for the centralizer of $N$ in $G$.
For $x\in\g$ and $s\in G$, we write $s{\cdot}x$ 
in place of $(\Ad s)x$. $K^o$ is the identity component of an algebraic group
$K$. If $\ah$ is a Lie algebra, then $\fe(\ah)\subset \ah\oplus\ah$ is the 
{\it commuting variety}, i.e., the set of all pairs of commuting elements.
We write just $\fe$ in place of $\fe(\g)$. Our general reference for 
nilpotent orbits is \cite{CoMc}.
\vno{.8}%
{\small {\bf Acknowledgements.} I am grateful to V.\,Ginzburg for kind 
information about his results and interesting discussions. Thanks are 
also due to A.\,Elashvili for friendly encouragement and to
 A.\,Broer for drawing my attention to \cite{nappes}.
This research was supported in part by RFFI Grant {\rus N0}\, 98--01--00598.}
\sekt{Principal nilpotent pairs\nopagebreak}%
We first review basic structure results on \pn pairs proved 
in \cite{vitya}.
\begin{rem}{Definition \ {\ququ (Ginzburg)}} \label{pnp}
A pair $\eb=(e_1,e_2) \in \g\times\g$
is called a {\it principal nilpotent pair} if the following holds:
\par
(i) $[e_1, e_2]=0$ \ and   $ \dim\z_\g(\eb)=\rk\g$;\par
(ii) For any $(t_1,t_2)\in \Bbbk^*\times \Bbbk^*$,  there exists 
$g=g(t_1,t_2)\in G$  such that \par \quad
$(\,t_1e_1,\,t_2e_2\,)=(g{\cdot}e_1,\,g{\cdot}e_2\,).$
\end{rem}%
The first step in Ginzburg's theory is that condition (ii) is equivalent to
the following one: there exists
an (associated semisimple) pair $\hb=(h_1, h_2)\in\g\times\g$ such that 
$\ad h_1$ and $\ad h_2$ have rational eigenvalues and
\beq [h_1, h_2]=0,\quad [h_i, e_j]=\delta_{ij} e_j \quad
(i,j\in \{1,2\})\ .  \label{comrel}
\eeq
In particular, the pair $\eb$ is nilpotent in the sense of introduction.
This $\hb$ determines the 
bi-grading of $\g$:
$\g_{k_1,k_2}=\{x\in\g\mid [h_j,x]=k_jx,\quad j=1,2\}$ and the induced
grading of $\z_\g(\eb)$. 
\begin{s}{Theorem {\ququ (see \cite[1.2]{vitya})}}   \label{smes}  \\
1. $\z_\g(\hb)$ is a Cartan subalgebra of $\g$; \\
2.  the eigenvalues of $\ad h_1$, $\ad h_2$ in $\g$ are integral; \\
3. $\displaystyle\z_\g(\eb)=\bigoplus_{i,j \in{\Bbb Z}_{\,\ge 0},
   (i,j)\neq(0,0)}\,\;\z_\g(\eb)_{i,j}$, i.e.,
$\z_\g(\eb)$ is graded by the `positive quadrant' without origin; \\
4. $\hb$ is determined uniquely up to conjugacy by $Z_G(\eb)^o$ (that is, the set 
of associated semisimple pairs form a single $Z_G(\eb)^o$-orbit). \qu
\end{s}%
Because of the last property it is natural to work with 
a (fixed) quadruple
$(\eb,\hb)$ rather than with the pair $\eb$.
Denoting $\el_i:=\z_\g(h_i)$ ($i=1,2$), 
we get $e_1,h_1\in\el_2$ and $e_2,h_2\in\el_1$.
Having the ${\Bbb Z}^2$-grading of $\g$ determined by $\hb$, 
one immediately sees 2 natural parabolic
subalgebras containing $\el_1$ and $\el_2$:\quad
$\p_1:=\displaystyle\bigoplus_{k_1\ge 0}\g_{k_1,k_2}=
\g_{\ge 0,\ast}$ (the right half-plane) and
$\p_2:=\displaystyle\bigoplus_{k_2\ge 0}\g_{k_1,k_2}=
\g_{\ast,\ge 0}$ (the upper half-plane). 
Then $\el_i$ is a Levi subalgebra of $\p_i$ and $e_i$ lies in the nilpotent
radical $(\p_i)^{nil}$ of $\p_i$. The main structure result is
\begin{s}{Theorem {\ququ (see \cite[sect.\,3]{vitya})}} \label{main1}
If $\eb$ is a \pn pair, then \\
{\rm (i)} $e_i$ is a Richardson element in $(\p_i)^{nil}$ (equivalently,
$\p_i$ is a polarization of $e_i$), $i=1,2$; \\
{\rm (ii)} $e_1$ (resp. $e_2$) is a regular nilpotent element in $\el_2$
(resp. $\el_1$).
\qu
\end{s}%
That the theory of \pn pairs has rich content follows from the classification
of such pairs in ${\frak sl}_N$, see \cite[5.6]{vitya}. In particular, the
following holds: given a nilpotent element $e\in{\frak sl}_N$, there exists
$e'$ such that $(e,e')$ is a \pn pair. The partition corresponding to $e'$
is conjugate to that for $e$. An explicit description of this pair is given
in terms of the corresponding Young diagram. This shows $\g$ may contain
many \pn pairs. Nevertheless, the following fundamental result is true:
\begin{s}{Theorem {\ququ (see \cite[3.9]{vitya})}} \label{main2}
The number of $G$-orbits of principal nilpotent pairs in $\g$ is finite.
\qu \end{s}%
Therefore the \pn pairs in simple Lie algebras can be classified.
(Actually, it may happen that $\g$ contains no non-trivial 
\pn pairs at all, see e.g.
$\GR{C}{2}$ or $\GR{B}{3}$ or $\GR{G}{2}$.) Some results 
to this end are found in \cite{wir}.

\sekt{Almost principal nilpotent pairs\nopagebreak}%
In this section we show that a large portion of Ginzburg's theory in the first
half of \cite{vitya}
remains valid in a more general setting. A motivation for this work came
from studying dual pairs associated with nilpotent pairs, see section 4.
Although some new phenomena arise, most of the proofs of \re{alm1} and
\re{alm3} are adapted from Ginzburg's ones.
\begin{rem}{Definition} \label{apnp}
A pair $\eb=(e_1,e_2) \in \g\times\g$
is called an {\it almost principal nilpotent pair} if the following holds:
\par
(i) $[e_1, e_2]=0$ \ and   $ \dim\z_\g(\eb)=\rk\g+1$;\par
(ii) there exists a pair of semisimple elements
 $\hb=(h_1, h_2)\in\g\times\g$ such that 
$[h_1, h_2]=0$ and $[h_i, e_j]=\delta_{ij} e_j \quad
(i,j\in \{1,2\})$.
\end{rem}%
Each pair $\hb$ satisfying condition (ii) is called an {\it associated 
semisimple\/} pair.
As in section 1, we shall consider the bi-grading
$\g=\oplus\g_{i,j}$ determined
by $\hb$.
For any subspace $M\subset\g$, one may define 3 filtrations:
\par 
$\bullet$ $e_1$-filtration: 
$M(i,\ast)=\{x\in M\mid (\ad e_1)^{i+1}x=0\}$, $i\ge 0$\ ;
\par 
$\bullet$ $e_2$-filtration: 
$M(\ast,j)=\{x\in M\mid (\ad e_2)^{j+1}x=0\}$, $j\ge 0$\ ;
\par 
$\bullet$  the double filtration: 
$M(i,j)=M(i,\ast)\cap M(\ast,j)$. \\
Following an idea of R.\,Brylinski, define the corresponding limits:\par
lim$_{e_1} M=\sum_{i\in {\Bbb Z}_{\ge 0}} (\ad e_1)^i M(i,\ast)\subset\g$, \par
lim$_{e_2} M=\sum_{j\in {\Bbb Z}_{\ge 0}} (\ad e_2)^j M(\ast,j)\subset\g$, \par
lim$_\eb M=\sum_{i,j\in {\Bbb Z}_{\ge 0}} (\ad e_1)^i(\ad e_2)^j M(i,j)\subset\g$.
\begin{s}{Theorem}  \label{alm1}
Let $\eb$ be an almost \pn pair and $\hb$ an associated semisimple pair. 
Then \par
{\rm (i)} $\z_\g(\hb)$ is a Cartan subalgebra of $\g$; \par
{\rm (ii)} $\z_\g(\eb)\cap\z_\g(\hb)=0$ or, equivalently, $\z_\g(\eb)_{0,0}=0$.
\end{s}%
\begin{proof*} 
We use an algebraized version of arguments in
\cite[sect.\,1]{vitya}. \\[.5ex]
(i) Consider the double filtration for $\te:=\z_g(\hb)$. 
Since $\te=\g_{0,0}$ and 
$(\ad e_1)^i(\ad e_2)^j\te(i,j)\subset \g_{i,j}$, the sum in the
definition is actually direct. Obviously, $\displaystyle 
\mbox{lim\,}_\eb\te\subset
\bigoplus_{i,j\in {\Bbb Z}_{\ge 0}}\z_\g(\eb)_{i,j}$. It follows from the
definition of double filtration
that $\te(i{-}1,j)\cap\te(i,j{-}1)=\te(i{-}1,j{-}1)$. Therefore
\[
\dim\,(\ad e_1)^i(\ad e_2)^j\te(i,j)=\dim\te(i,j)-\dim\te(i{-}1,j)-
\dim\te(i,j{-}1)+\dim\te(i{-}1,j{-}1) \, .
\]
Now, an easy summation proves that $\dim\,(\mbox{lim\,}_\eb\te)=\dim\te$. Thus,
\[
\rk\g\le\dim\te\le\dim\,\bigl(
\bigoplus_{i,j\in {\Bbb Z}_{\ge 0}}\z_\g(\eb)_{i,j}
\bigr) \le\dim\z_\g(\eb)=\rk\g+1\, .
\]
Since $\te$ is a Levi subalgebra, $\dim\te-\rk\g$ is 
even. Hence $\te$ must be a Cartan subalgebra. \\[.5ex]
(ii) Assume that $h$ is a nonzero element in $\z_\g(\eb)\cap \te$.
Then $e_1,e_2$ lie in the Levi subalgebra $\el:=\z_\g(h)$.
By a result of Richardson \cite{roger}, $\fe(\el)$ is irreducible and
the pairs of semisimple elements are dense in $\fe(\el)$. Therefore
\[
\rk\g=\rk\el\le\dim\z_\el(e_1,e_2)\le\dim\z_\g(e_1,e_2)=\rk\g+1\ .
\]
Associated with $\el$, there is a decomposition $\g=\n_+\oplus\el\oplus\n_-$, 
where $[\el,\n_\pm]=\n_\pm$. 
It follows that 
$\z_\g(\eb)=\z_{\n_-}(\eb)\oplus\z_\el(\eb)\oplus\z_{\n_+}(\eb)$ and
$\dim\z_\g(\eb)=2\dim\z_{\n_+}(\eb)+\dim\z_\el(\eb)$. Obviously, the first
summand is positive and we obtain $\dim\z_\g(\eb)\ge\rk\g+2$.
This contradiction proves the claim (ii).
\end{proof*}%
\begin{s}{Corollary}  \label{alm2} We have
${\lim\,}_\eb\te=\displaystyle
\bigoplus_{i,j\in {\Bbb Z}_{\ge 0}, (i,j)\ne (0,0)}\z_\g(\eb)_{i,j}$. 
In particular,\\
$(\ad e_1)^i(\ad e_2)^j\te(i,j)=\z_\g(\eb)_{i,j}$ for all 
$i,j\in {\Bbb Z}_{\ge 0}$.
\end{s}\begin{proof}
It is already proved that the inclusion ``$\subset$" holds. 
Since $\te$ is Cartan and the pair $\eb$ is not principal, 
it follows from \cite[1.13]{vitya} that 
$\z_\g(\eb)\ne\bigoplus_{i,j\in {\Bbb Z}_{\ge 0}, (i,j)\ne (0,0)}
\z_\g(\eb)_{i,j}$.
Then the assertion follows for dimension reason. 
\end{proof}%
Unlike the case of \pn pairs (see  \re{smes}), the eigenvalues of 
$\ad h_1$ and $\ad h_2$ are not necessarily integral and $\z_\g(\eb)$
is not necessarily graded by `positive quadrant'. As we shall see in 
\re{alm4}, these two
conditions form a dichotomy in case of almost \pn pairs.
\begin{s}{Theorem}  \label{alm3} 
Let $\eb$ be an (almost) \pn pair with an associated semisimple pair
$\hb$. Put $\el_i=\z_\g(h_i)$ and let $\ce_i$ denote the centre of $\el_i$
($i=1,2$). Then \par
1. $e_1$ is a \rn element in $\el_2$ and 
 $e_2$ is a \rn element in $\el_1$; \par
2. ${\lim\,}_{e_1}\te=\z_\g(e_1,h_2)$ and 
${\lim\,}_{e_2}\te=\z_\g(e_2,h_1)$;\par
3. $\dim\z_\g(e_1,h_1,e_2)=\dim\z_\g(e_1,h_1,h_2)$ and 
$\dim\z_\g(e_2,h_2,e_1)=\dim\z_\g(e_2,h_2,h_1)$; \par
4. $\z_\g(e_1,h_1,h_2)=\ce_2$ and $\z_\g(e_2,h_1,h_2)=\ce_1$.
\end{s}\begin{proof}
By symmetry, it suffices to prove the first half of each item.
The proof applies to both \pn\ and almost \pn pairs.\\
1 \& 2. This part is essentially the same as in \cite{vitya}. 
Consider the $e_1$-limit:
${\mathrm {lim}\,}_{e_1}\te=\sum_{i\ge 0} (\ad e_1)^i\te(i,\ast)$,
which lies in $\z_\g(e_1,h_2)$. Since different summands have different weights
relative to $\ad h_1$, the sum is direct 
and therefore
$\rk\g\le\dim\z_\g(e_1,h_2)$. The space $\z_\g(e_1,h_2)$ possesses 
the $e_2$-filtration and $\mbox{lim\,}_{e_2}\z_\g(e_1,h_2)\subset\z_\g(\eb)$.
For similar reason, $\dim\,(\mbox{lim\,}_{e_2}\z_\g(e_1,h_2))=\z_\g(e_1,h_2)$ and
hence
$\dim\z_\g(e_1,h_2)\le\rk\g+1$. As in the proof of \re{alm1}(i), one may 
conclude by making use of 
the parity argument: $e_1$ lies in the reductive Lie algebra $\el_2$ and
therefore $\dim\z_\g(e_1,h_2)=\dim\z_{\el_2}(e_1)$ must have the same parity
as $\rk\g=\rk\el_2$. \\
3. Applying the formula in \re{alm2} with $i=0$ gives
\[
\bigoplus_j (\ad e_2)^j\te(0,j)=\bigoplus_j\z_\g(\eb)_{0,j}=
\z_\g(e_1,e_2,h_1) \ .
\]
Obviously, dimension of the left-hand side is 
$\sum_j \dim\,(\te(0,j)/\te(0,j{-}1))=
\dim\te(0,\ast)$. Since $\te(0,\ast)=\z_\g(e_1,h_1,h_2)$, we are done.\\
4. 
Since $e_1,h_1\in\z_\g(h_2)=\el_2$,
we have $\z_\g(e_1,h_1,h_2)\supset\ce_2$. By either \re{main1}(ii) or 
\re{alm3}(1), 
$e_1$ is a \rn element in
$\el_2$. Therefore $\z_{\el_2}(e_1)= \z_\g(e_1,h_2)=\ce_2\oplus\n$, 
where $\n\subset [\el_2,\el_2]$ consists of nilpotent elements.
Finally, $\z_\g(e_1,h_1,h_2)\subset
\z_\g(h_1,h_2)$ and therefore $\z_\g(e_1,h_1,h_2)$ consists of 
semisimple elements. Whence $\z_\g(e_1,h_1,h_2)=\ce_2$.
\end{proof}%
\indent By \re{alm2}, 
$\z_+:=\displaystyle
\bigoplus_{i,j\in {\Bbb Z}_{\ge 0},(i,j)\ne (0,0)}\z_\g(\eb)_{i,j}$
is of codimension one in $\z_\g(\eb)$, if $\eb$ is an almost \pn pair. 
Hence there is an ``extra" vector
$x$ in some $\g_{p,q}$ such that $\z_\g(\eb)=
\z_+\oplus\langle x\rangle$. 
We already know that $(p,q)\not\in ({\Bbb Z}_{\ge o})^2$. It follows from
\re{alm3}(2) that the eigenvalues of $\ad h_1$ (resp. $\ad h_2$) in
$\z_\g(e_1,h_2)$ (resp. $\z_\g(e_2,h_1)$\,) are nonnegative integers.
Therefore $x\not\in\el_i$ ($i=1,2$). That is, $pq\ne 0$.
\begin{s}{Theorem}   \label{alm4}
1. There are 2 mutually exclusive possibilities for $p,q$. Either \\[.5ex]
\hspace*{1.1cm} $({\Bbb Z})$\hspace{1cm} $p,q\in \Bbb Z$ and $pq<0$,\quad
or \\
\hspace*{.9cm}$(\mbox{non-}{\Bbb Z})$\hspace{.5cm}
$p,q\in \frac{1}{2}{\Bbb Z}\setminus\Bbb Z$ and
$p,q>0$. \\
2.  In both cases, 
$\z_\g(\eb)$ is nilpotent and contains no semisimple elements. Moreover,
$\z_\g(\eb)$ is Abelian in the non-$\Bbb Z$ case.
\end{s}\begin{proof*}
1. For ($\Bbb Z$): Suppose that $p,q\in\Bbb Z$, i.e., all the eigenvalues
of $\hb$ in $\z_\g(\eb)$ are integral.  We need to prove here that the
case $p<0,\ q<0$ is impossible. Assume not and $p_0=-p>0$,
$q_0=-q>0$. A standard calculation with the Killing form on $\g$ shows that
$\z_\g(\eb)_{p,q}\ne 0$ if and only if
$\g_{p_0,q_0}\not\subset \mbox{Im\,}(\ad e_1)+\mbox{Im\,}(\ad e_2)$.
By definition, put ${\cal D}=
\g_{p_0,q_0}\setminus (\mbox{Im\,}(\ad e_1)+\mbox{Im\,}(\ad e_2))$. For each
$y\in\cal D$, consider the finite set $I_y=\{(k,l)\in ({\Bbb Z}_{\ge 0})^2\mid
(\ad e_1)^k(\ad e_2)^ly\ne 0\}$, with the lexicographic ordering.
This means $(k,l)\prec(k',l')$ $\Leftrightarrow$ 
$k<k'$ or $k=k'$ and $l<l'$. Denote by $m(I_y)$ the unique maximal element
in $I_y$. Let $y^*\in\cal D$ be an element such that
$(k_0,l_0):=m(I_{y^*})\preceq m(I_z)$ for all $z\in\cal D$. Then
$(\ad e_1)^{k_0}(\ad e_2)^{l_0}y^*$ is a nonzero element in $\z_\g(\eb)\cap
\g_{p_0+k_0,q_0+l_0}$. By \re{alm2}, there is $t\in\te(p_0+k_0,q_0+l_0)$
such that $(\ad e_1)^{p_0+k_0}(\ad e_2)^{q_0+l_0}t=
(\ad e_1)^{k_0}(\ad e_2)^{l_0}y^*$. Then $(\ad e_1)^{k_0}(\ad e_2)^{l_0}(y^*-
(\ad e_1)^{p_0}(\ad e_2)^{q_0}t)=0$. Since $p_0>0,q_0>0$, we have
$z^*=y^*-(\ad e_1)^{p_0}(\ad e_2)^{q_0}t$ is nonzero and belongs to $\cal D$.
However, $I_{z^*}\subset I_{y^*}\setminus\{(k_0,l_0)\}$. Therefore
$m(I_{z^*})< m(I_{y^*})$, which contradicts the choice of $y^*$.
Thus, the case $p<0,q<0 $ is impossible. \par
For (non-$\Bbb Z$): Suppose $(p,q)\not\in {\Bbb Z}\oplus\Bbb Z$. Consider
the set ${\cal J}=\{(k,l)\mid \g_{k,l}\ne 0\mbox{ and } (k,l)
\not\in {\Bbb Z}\oplus{\Bbb Z}\}$. Because $\langle x\rangle$ is the unique 
``non-integral" homogeneous subspace of $\z_\g(\eb)$, ${\cal J}$ lies in the
single coset space $(p,q)+({\Bbb Z}\oplus\Bbb Z)$ and has a unique 
``north-east" corner. Obviously, $(p,q)$ is this corner. Since $\dim\g_{m,n}=
\dim\g_{-m,-n}$ for all $(m,n)$, this corner must lie in the positive
quadrant. The condition $(-p,-q)\in (p,q)+({\Bbb Z}\oplus\Bbb Z)$ implies
$p,q\in \frac{1}{2}\Bbb Z$. It remains to demonstrate that {\it both\/} 
$p,q$ must be fractional.  Assume not, and $p\in\Bbb Z$, while $q$ is
fractional. Consider a
``path inside of $\g$" connecting the points $(-p,-q)$ and $(p,q)$:
Starting from a nonzero element in $\g_{-p,-q}$, we may always apply either
$\ad e_1$ or $\ad e_2$ until we arrive at $\ap x\in\g_{p,q}$ ($\ap\ne 0$). 
Since $p$ is integral, 
we must intersect somewhere the vertical axis. This means $\ad h_2$ has
a fractional eigenvalue in $\el_1$. It then follows from nilpotency of
$\ad e_2$ that $\ad h_2$ has a fractional eigenvalue in $\z_{\el_1}(e_2)$
as well. However, this contradicts to \re{alm3}(2).
\\
2. The pairs $(k,l)$ such that $\z_\g(\eb)_{k,l}\ne 0$ are said to be
{\it bi-weights\/} of $\z_\g(\eb)$. 
In either case, the bi-weights lie in an open half-plane of $\Bbb Q\oplus
\Bbb Q$, hence the assertion. In the non-$\Bbb Z$ case, $(p,q)$ is the unique 
nonintegral bi-weight. Since $(0,0)$ is not a bi-weight (see \ref{alm1}(ii)),
this implies $[\z_+,\, x]=0$.
It is also easily seen that $\z_+=\mbox{lim\,}_\eb\te$ is Abelian.
\end{proof*}%
\begin{s}{Corollary}  \label{half-int}
If $\hb$ is any associated semisimple pair then the eigenvalues of $\ad h_1$,
$\ad h_2$ in $\g$ are at least half integers. \qu
\end{s}%
An almost \pn pair is said
to be either of {\it $\Bbb Z$-type\/} or of {\it non-$\Bbb Z$-type\/}
according to the two possibilities in Theorem
\ref{alm4}(1). It will be proved below that all associated semisimple pairs
are $Z_G(\eb)^o$-conjugate. Therefore the property of being of 
$\Bbb Z$-type 
does not depend on the choice of $\hb$. 
\begin{s}{Corollary}  \label{*inv}
Let $\eb$ be an almost \pn pair of non-$\Bbb Z$-type. Then there is an
inner
involution $\theta\in\mbox{Aut\,}\g$ such that $\g^\theta$ is semisimple and
$\eb$ is a \pn pair in $\g^\theta$.
\end{s}\begin{proof}
Define $\theta\in\mbox{Hom\,}(\g,\g)$ by 
$ \theta\mid_{\g_{i,j}}=\left\{ 
\begin{array}{rl} \mbox{id} & \mbox{ if $i,j\in\Bbb Z$}\\
-\mbox{id} & \mbox{ if $i,j\in {\Bbb Z}+\frac{1}{2}$} \end{array} \right.
$.
It is an inner automorphism of $\g$. Then $e_1,e_2\in\g^\theta$, 
$\rk\g^\theta=\rk\g$, and $\dim\z_{\g^\theta}(\eb)=\rk\g^\theta$. As $\z_{\g^\theta}(\eb)$
contains no semisimple elements, $\g^\theta$ is semisimple.
\end{proof}%
It is worth noting that the two cases in \re{alm4}  really occur:
\begin{rem}{Example}  \label{2-vozm}
Take $\g={\frak sp}_4$. Let $\ap=\esi_1-\esi_2$ and
$\beta=2\esi_2$ be the usual simple roots. Denote by
$e_\mu$ a nonzero root vector corresponding to $\mu$. Then
$(e_{2\ap+\beta},e_\beta)$ is an almost \pn pair of $\Bbb Z$-type and
$(e_{\ap+\beta},e_{2\ap+\beta})$ is an almost \pn pair of 
non-$\Bbb Z$-type.
In both cases, $\z_\g(\eb)=\langle e_{2\ap+\beta},e_{\ap+\beta},e_\beta
\rangle$, but associated semisimple pairs  are essentially different.
The details are left to the reader.
\end{rem}%
As in section 1, define the parabolic subalgebras $\p_1$ and $\p_2$. Unlike
the \pn case, $e_i$ is not necessarily a Richardson element in the nilpotent
radical $(\p_i)^{nil}$ of $\p_i$. The precise statement is as follows.
\begin{s}{Theorem}  \label{alm-rich}
Let $\eb$ be an almost \pn pair. \\
{\rm (i)} Suppose $\eb$ is of non-$\Bbb Z$-type. Then neither of $e_i$'s is
Richardson in $(\p_i)^{nil}$; \\
{\rm (ii)} Suppose $\eb$ is of $\Bbb Z$-type with, say, $q>0$ and $p<0$. Then
$e_2$ is Richardson in $(\p_2)^{nil}$, while $e_1$ is not Richardson in
$(\p_1)^{nil}$.
\end{s}\begin{proof}
(i) Let $(k_0,l_0)$ be the minimal element in
$\bigl( (p,q)+({\Bbb Z}\oplus {\Bbb Z})\bigr)\cap({\Bbb Q}_{>0}\oplus
{\Bbb Q}_{>0})$ with respect to the lexicographic ordering. Then
$\g_{k_0,l_0}\subset (\p_i)^{nil}$, while $\g_{k_0,l_0}\not
\subset [\p_i,e_i]$ ($i=1,2$). It is not hard to prove that $(k_0,l_0)=
(1/2,1/2)$, but we do not need this. \par
(ii) Now the eigenvalues of $\ad \hb$ are integral and the bi-weights
of $\z_\g(\eb)$ lie in the upper half-plane. The same argument as in 
\cite[1.12]{vitya} shows that $\ad e_2 : \g_{\ap,\beta}\rar
\g_{\ap,\beta+1}$ is injective for all $\ap$ and $\beta<0$. (Otherwise
we would find an element $0\ne y\in\z_\g(\eb)_{\nu,\beta}$ with $\nu\ge\ap$, 
$\beta<0$.) Then, by duality, $\ad e_2$ is surjective 
for $\beta\ge 0$. In particular, $[\p_2,e_2]=[\g_{\ast,\ge 0},e_2]=
\g_{\ast,\ge 1}=(\p_2)^{nil}$. \\
On the other hand, $\ad e_1: \g_{p,q}\rar \g_{p+1,q}$ is not injective.
Hence $\ad e_1: \g_{-p-1,-q}\rar\g_{-p,-q}$ is not surjective, i.e., 
$[\p_1,e_1]=[\g_{\ge 0,\ast},e_1]\ne\g_{\ge 1,\ast}=(\p_1)^{nil}$.
\end{proof}%
Recall the notion, due to Lusztig and Spaltenstein, of a {\it special\/} 
nilpotent 
orbit. Let $\N/G$ be the set of all nilpotent orbits in $\g$. The closure
ordering ``$\co_1\cyeq \co_2\ \Leftrightarrow\ \co_1\subset\ov{\co_2}$"
makes $\N/G$ a finite poset. In \cite[ch.\,III]{Sp}, 
N.\,Spaltenstein studied a 
duality in $(\N/G,{\cyeq})$. He proved that there exists an order-reversing
mapping $d:\N/G\rar \N/G$ such that \par
 (a) $\co\cyeq d^2(\co)$ for all $\co\in\N/G$; \par
 (b) For any Levi subalgebra $\el\subset\g$, $d$ takes the $G$-orbit through 
the regular nilpotent elements in $\el$ to the Richardson orbit associated to
$\el$.
\\[.5ex]
Such a mapping can uniquely be determined, in a purely combinatorial way,
for the classical Lie algebras and for $\GR{E}{7}$. In the remaining cases,
a natural choice among finitely many possibilities can be done. Then one 
of the definitions of specialness is that $(\N/G)_s:=d(\N/G)$ is just the
set of special orbits. An important feature of $(\N/G)_s$ is that
$d\mid_{(\N/G)_s}$ is an order-reversing involution. In case of ${\frak sl}_n$,
this is just the conjugation on the set of all partitions of $n$. With these
results at hand, an immediate consequence of the previous theorem is:
\begin{s}{Proposition}  \label{non-special}
Let $\eb$ be an almost \pn pair of $\Bbb Z$-type, as in \re{alm-rich}(ii).
Then $Ge_1$ is not special.
\end{s}\begin{proof*}
In view of Theorem \ref{alm3}(1), the assertion \re{alm-rich}(ii) can be
restated as: $d(Ge_1)=Ge_2$ and $d(Ge_2)\ne Ge_1$. Assume now that
$Ge_1$ is special, i.e. $Ge_1=d(\co)$ for some $\co\in\N/G$. 
Then $d^2(\co)=Ge_2$ and 
$Ge_1=d(\co)=d^3(\co)=d(Ge_2)$. A contradiction! 
\end{proof*}
\begin{s}{Corollary}  \label{sl-special}
There are no almost \pn pairs in ${\frak sl}_n$.
\end{s}\begin{proof}
1. By Corollary \ref{*inv}, any almost \pn pair of non-$\Bbb Z$-type
yields an inner involution $\theta$ such that $\g^{\theta}$ is semisimple.
But ${\frak sl}_n$ has no such involutions. \\
2. Since all nilpotent orbits in ${\frak sl}_n$ are Richardson and
hence special, there are no almost \pn pairs of $\Bbb Z$-type as well.
\end{proof}%
The following easy result is needed in the proof of \re{conj}.
\begin{s}{Lemma}  \label{lem}
Let $h_1,h_2$ be two commuting semisimple elements. Let $\n\subset\g$ be
a subspace such that $[h_i,\n]\subset\n$ $(i=1,2)$ and 
$\n\cap\z_\g(h_1,h_2)=\{0\}$.
Then $\dim\{(n_1,n_2)\subset\n\oplus\n\mid [h_1,n_2]=[h_2,n_1]\}=\dim\n$.
\qu \end{s}%
\vskip-1.5ex
\begin{s}{Theorem}   \label{conj}
Let $\eb$ be an almost \pn pair. Let $\hb$ and $\hb'=(h_1',h_2')$ be two 
associated semisimple pairs. Then there exists $u\in Z_G(\eb)^o$ such that
$u{\cdot}h_i=h_i'$ ($i=1,2$).
\end{s}\begin{proof}
Let $\cw$ be the set of all associated semisimple pairs.
Obviously, $Z_G(\eb)^o{\cdot}\hb\subset\cw$. 
It follows from \re{alm4} that $Z_G(\eb)^o$ is unipotent and therefore
$Z_G(\eb)^o{\cdot}\hb$ is closed in $\g\oplus\g$.
Since $Z_G(\eb)\cap Z_G(\hb)$ is finite, $\dim Z_G(\eb)^o{\cdot}\hb=\rk\g+1$. 
On the other hand,
$h_i'-h_i\in\z_\g(\eb)$ ($i=1,2)$. Therefore $\cw\subset
(h_1+\z_\g(\eb),h_2+\z_\g(\eb))\cap\fe=:\tilde\cw$. Recall that $\fe$ is the
commuting variety. Hence, the assertion is equivalent to that $\tilde\cw$ is 
irreducible and of dimension $\rk\g+1$.
Our next analysis relies on the structure of $\z_\g(\eb)$ described in
\re{alm2},\re{alm4}. We have $\z_\g(\eb)=\z_+\oplus\langle x\rangle$, 
where $\z_+=\mbox{lim\,}_\eb\te$
and $x\in\g_{p,q}$. In both cases described in \re{alm4}, one has
$\z_+$ is Abelian and
$[x,\z_+]\subset\z_+$. 
Let $(h_1+n_1+\nu x,h_2+n_2+\tau x)\in\tilde\cw$, where $n_i\in\z_+$, 
$\nu,\tau\in\Bbbk$. The $x$-coordinate of the commutator is equal to
$\tau p-\nu q$. Therefore $(\nu,\tau)=c(p,q)$ with $c\in\Bbbk$.
Vanishing of the $\z_+$-component yields the equation
\[
[h_1,n_2]-[h_2,n_1]+c(q[n_1,x]-p[n_2,x])=0 \ .
\]
Having fixed $c$, we obtain a system of linear equations for $n_1,n_2$. More
precisely, consider the family of linear mappings \par
$\nu_c: \z_+\oplus\z_+\rar\z_+$,
$(n_1,n_2)\mapsto [h_1,n_2]-[h_2,n_1]+c(q[n_1,x]-p[n_2,x])$. \\
Then
$\tilde\cw=\sqcup_{c\in\Bbbk}(\hb+\mbox{Ker\,}\nu_c+c(px,qx))$. 
By Lemma \ref{lem},
$\mbox{Ker\,}\nu_0\simeq\{(n_1,n_2)\mid [h_1,n_2]=[h_2,n_1]\}$ is of
dimension $\dim\z_+=\rk\g$. That is, $\nu_0$ is onto. It follows that
$\dim\mbox{Ker\,}\nu_c=\rk\g$ for all but finitely many $c\in\Bbbk$.
Therefore $\tilde\cw$ has a unique irreducible component passing through
$(h_1,h_2)$ and $\mbox{dim\,}_\hb\tilde\cw=\rk\g+1$. Recall that the bi-weights
of $\z_\g(\eb)$ lie in an open half-space in $\Bbb Q\oplus\Bbb Q$.
Therefore there exists a 1-parameter subgroup in the maximal torus
$Z_G(\hb)$ which contracts everything in the affine space
$(h_1+\z_\g(\eb),h_2+\z_\g(\eb))\subset\g\oplus\g$ to the point $\hb$. 
Hence $\tilde\cw$ is
a cone with vertex $\hb$. Thus, $\tilde\cw$ is irreducible and of dimension
$\rk\g+1$.
\end{proof}%
While $Z_G(\eb)$ is always connected in case of \pn pairs 
(see \cite[3.6]{vitya}), connectedness in the almost principal case depends
on the type.
\begin{s}{Proposition}  \label{con-apn}
Let $\eb$ be an almost \pn pair. Then \par
1. $Z_G(\eb)$ is connected, if $\eb$ is of $\Bbb Z$-type; \par
2. $Z_G(\eb)$ is disconnected, if $\eb$ is of non-$\Bbb Z$-type.
\end{s}\begin{proof*}  From Theorem \ref{alm4}(2), it follows that 
$Z_G(\eb)$ is a semi-direct product
of the unipotent group $Z_G(\eb)^o$ and a finite group $F$.  \\
1. Take an
arbitrary $s\in F$. It is a semisimple element of finite order. Since
$s{\cdot}\hb$ is an associated semisimple pair for $\eb$, it follows from
\re{conj} that $s{\cdot}\hb=u{\cdot}\hb$ for some $u\in Z_G(\eb)^o$.
Hence $t:=s^{-1}u\in Z_G(\hb)=T$. By \ref{alm-rich}(ii), one may assume that
$e_2$ is Richardson in $(\p_2)^{nil}$. Since $t{\cdot}e_1=e_1$ and $e_1$ is
\rn in $\el_2$ (see \ref{alm3}), $t$ is in the centre of $Z_G(h_2)=:L_2$. 
Because $\el_2$ and $e_2$ generate the parabolic subalgebra $\p_2$ and
$t{\cdot}e_2=e_2$, we get 
$t{\cdot}z=z$ for any $z\in\p_2$. This clearly implies that
$t$ is in the centre of $G$. Since $G$ is adjoint, we obtain
$s=u=1\in G$. \\
2. By Corollary \ref{*inv}, $Z_G(\eb)$ contains a semisimple element of 
order two.
\end{proof*}%
\begin{rem}{Example}  \label{seriya}
Here we describe a series of almost \pn pairs in
symplectic Lie algebras. Let $\g={\frak sp}_{4n}={\frak sp}(\V)$ and let
$v_1,\dots,v_{4n}$ be a basis of $\V$ such that the $\g$-invariant 
skew-symmetric form is $B(z,y)=z_1y_{4n}+\dots+z_{2n}y_{2n+1}-
z_{2n+1}y_{2n}-\dots-z_{4n}y_{1}$. Define the operators $e_1,e_2\in{\frak sp}(\V)$
by the formulas: \par
$e_1(v_j)=v_{j-2}$ ($j\ge 2n+1$), $e_1(v_j)=-v_{j-2}$ ($3\le j\le 2n$);\par
$e_2(v_{2j})=v_{2j-3}$ ($j\ge n+1$), $e_2(v_{2j})=-v_{2j-3}$ 
($2\le j\le n$).\\
If $e_i(v_j)$ is not specified, this means it is equal to zero.
The orbit $G{\cdot}e_1$ (resp. $G{\cdot}e_2$) corresponds to the 
partition $(2n,2n)$ 
(resp. $(2,\dots,2,1,1)$).
Then $[e_1,e_2]=0$ and $\z_\g(e_1,e_2)=\langle
e_1,\,e_1^3,\dots,e_1^{2n-1},\,e_2,\,e_1^2e_2,\dots,e_1^{2n-2}e_2,x \rangle
$, where $x$ is the operator taking $v_{4n-1}$ to $v_2$.
Hence $\z_\g(e_1,e_2)$ is Abelian and its dimension is $2n+1$. 
An associated semisimple pair consists of
$h_1=\mbox{diag\,}(t_1,\dots,t_{4n})$, where $t_{2i}=n{+}1{-}i$,
$t_{2i-1}=n-i$ ($i=1,\dots,2n$), and 
$h_2=\mbox{diag\,}(1/2,-1/2,1/2,-1/2,\dots)$. The bi-weights of 
$\z_\g(e_1,e_2)$ are: \\
\hspace*{1cm}
$(1,0),\,(3,0),\dots,(2n{-}1,0),\,(0,1),\,(2,1),\dots,(2n{-}2,1),\,(2n,-1)$,\\
where the ordering corresponds to that of basis vectors. Therefore these 
almost \pn pairs are of $\Bbb Z$-type.
Note that for $n=1$ we obtain one of the pairs given in \re{2-vozm}.
See also example in \re{nonZ}.
\end{rem}%
{\bf Remarks}. Here we collect some observations, which require either a 
longer presentation or a further development. \par
{\sf 1.} It is true that the number of $G$-orbits of almost \pn pairs is
finite (cf. \ref{main2}). This will be published elsewhere. \par
{\sf 2.} Any simple Lie algebra contains a non-trivial nilpotent pair
$\eb$ such that
$\dim\z_\g(\eb)\le \rk\g+1$.
\par
{\sf 3.} 
All known examples of almost \pn pairs occur in 
$\GR{B}{m}$, $\GR{C}{m}$, $\GR{G}{2}$. Apparently, this means that 
almost \pn pairs do not exist in the simply-laced case. \par
{\sf 4.}
It is likely that $\z_\g(\eb)$ is also Abelian for the almost \pn pairs
of $\Bbb Z$-type, but I did not succeed in finding a proof. 

\sekt{Rectangular nilpotent pairs\nopagebreak}%
A description of \pn pairs in ${\frak sl}_N$ obtained by
Ginzburg shows that in general $h_i\not\in\mbox{Im\,}(\ad e_i)$, i.e.,
 $\{e_i,h_i\}$ can not be included in a simple 3-dimensional
subalgebra. However, the theory becomes much simpler, if
this can be done, see \cite{wir}. This motivates the
following
\begin{rem}{Definition}  \label{rectangle}
A pair of nilpotent elements $(e_1,e_2)$ is called {\it rectangular\/} 
whenever there exists an $\tri$-triple, containing $e_1$,
that commutes with $e_2$.
\end{rem}%
Recall that an $\tri$-triple $\{e,\tih,f\}$ satisfies the commutator
relations $[\tih,e]=2e$, $[\tih,f]=-2f$, $[e,f]=\tih$. The famous 
Dynkin--Kostant theory describes conjugacy classes of $\tri$-triples and
the structure of $\z_\g(e)$ through the use
of $\{e,\tih,f\}$. (See either the original papers \cite{EBD} and \cite{ko59}
or a modern presentation in \cite[ch.\,4]{CoMc}.) 
Here are some results of this theory together with
related notions. 
The semisimple element $\tih$ is called a {\it characteristic} of $e$.
Given $e$ and $\tih$, the third member of $\tri$-triple is uniquely 
determined and $\z_\g(e,\tih)=\z_\g(e,\tih,f)$.
Let  $
\g=\bigoplus_{i\in {\Bbb Z}}\g(i)$ be the $\Bbb Z$-grading
defined by $\ad \tih$. Note that $\g(0)$ is nothing but $\z_\g(\tih)$. Then
$\z_\g(e)=\bigoplus_{i\ge 0}\z_\g(e)_i$ and $\z_\g(e)_0$ is a maximal 
reductive subalgebra in $\z_\g(e)$. Moreover, $\z_\g(e)_0=\z_\g(e,\tih,f)$. 
Putting $\z_\g(e)_{odd}=\bigoplus_{i\ odd}\z_\g(e)_i$ and likewise for 
``even", we have 
\beq  \label{izomy}
\z_\g(e)_{even}\simeq\g(0) \mbox{\quad and\quad }
\z_\g(e)_{odd}\simeq\g(1) \mbox{ as }\z_\g(e)_0\mbox{-module}.
\eeq
The element $e$ is called {\it even} whenever all the eigenvalues of
$\ad\tih$ are even. Obviously, $e$ is even if and only if $\g(1)=0$ if and only
if $\dim\z_\g(e)=\dim\z_\g(\tih)$. Then the weighted Dynkin
diagram of $e$ contains only numbers 0 and 2. An $\tri$-triple containing
regular elements is called {\it principal}; $e$ is regular if
and only if it is even and $\z_g(\tih)$ is a Cartan subalgebra. \\
Since all $\tri$-triples containing $e$ are $Z_G(e)^o$-conjugate, the above
properties have intrinsic nature.
\begin{s}{Lemma}  \label{triple} \\
1. The following conditions are equivalent for a  pair $\eb$
of nilpotent elements: \par
{\rm (i)} $(e_1,e_2)$ is rectangular; \par
{\rm (ii)} there exist commuting $\tri$-triples
$\{e_1,\tih_1,f_1\}$ and $\{e_2,\tih_2,f_2\}$.  \\
2. If $\eb$ is rectangular and $\g=\oplus\g_{ij}$ is the ${\Bbb Z}^2$-grading
defined by $(\tih_1,\tih_2)$, then $\z_\g(\eb)$ is graded by `positive 
quadrant'. \\
3. If $\eb$ is a rectangular (almost) \pn pair, then we may assume that 
$\hb=(\tih_1/2,\tih_2/2)$.
\end{s}%
\begin{proof} 1. 
Suppose an $\tri$-triple $\{e_1,\tih_1,f_1\}$ commutes with $e_2$.
Then we may choose an $\tri$-triple containing $e_2$ inside of the reductive
algebra $\z_\g(e_1,\tih_1,f_1)$.\\
2. This readily follows from the Dynkin-Kostant theory. \\
3. In this case $\tih_1/2,\tih_2/2$ satisfy
commutator relations \re{comrel}. From \re{smes}(4) and \re{conj}, we then
conclude that  $(\tih_1/2,\tih_2/2)$ is $Z_G(\eb)^o$-conjugate to $\hb$. 
\end{proof}%
Obviously, any rectangular pair is nilpotent in the sense of Introduction.
Because one may use the powerful ${\frak sl}_2$-machinery in the rectangular
case, it seems likely that any reasonable question concerning rectangular
pairs has an immediate answer. For instance, the following is proved in
\cite[Th.\,7.1]{wir}:
\begin{s}{Theorem}  \label{th-1A}
Let $\{e,\tih,f\}$ be an $\tri$-triple. Then $e$ is a member of a
rectangular \pn pair if and only if $e$ is 
even and a (any) \rn element in $\z_\g(e,\tih,f)$ is regular in
$\z_\g(\tih)$ as well. \qu
\end{s}%
It is not hard to find a similar statement in the almost principal case:
\begin{s}{Theorem} \label{rect-apnp}
Let $\{e,\tih,f\}$ be an $\tri$-triple. Then $e$ is a member of a
rectangular almost \pn pair if and only if the following holds: \par
1. a (any) \rn element in $\ka:=\z_\g(e,\tih,f)$ is also regular in
$\z_\g(\tih)$; \par
2. $e$ is \underline{not} even (i.e., $\g(1)\ne 0$) and $\dim\z_{\g(1)}(e')=1$,
if $e'\in\ka$ is regular nilpotent. \\
Under these hypotheses, if $e'\in\ka$ is regular nilpotent, then $(e,e')$ is
an almost \pn pair.
\end{s}\begin{proof}
The proof is much the same as for the previous assertion. Take a nilpotent
element $e'\in\ka$. It then follows from \re{izomy} that
\[ \dim\z_\g(e,e')=\dim\z_{\g(0)}(e')+\dim\z_{\g(1)}(e') \ .
\]
If $\dim\z_\g(e,e')=\rk\g+1$, then we must have
$\dim\z_{\g(0)}(e')=\rk\g(0)=\rk\g$ and $\dim\z_{\g(1)}(e')=1$ for 
``parity" reason. Thus $e'$ is regular in $\g(0)$ and hence in $\ka$.
This argument can be reversed.
\end{proof}%
Note that a rectangular almost \pn pair is necessarily of non-$\Bbb Z$-type 
and that
condition 2 can be restated as follows: $\g(1)$ is a simple
$\langle e',\tih',f'\rangle$-module, if 
$\{e',\tih',f'\}\subset\ka$ is a
principal $\tri$-triple.
\begin{rem}{Example}  \label{nonZ}
Let $\g={\frak sp}_{2n}$. For $0<k<n$,
consider the subalgebra ${\frak sp}_{2k}\oplus
{\frak sp}_{2n-2k}\subset {\frak sp}_{2n}$. 
Let $e_1$ (resp. $e_2$) be
a regular nilpotent element in ${\frak sp}_{2k}$ (resp. ${\frak sp}_{2n-2k}$).
Then $(e_1,e_2)$ is a rectangular almost \pn pair.
\end{rem}%
We refer to \cite{wir2} for a classification of rectangular 
almost \pn pairs in the simple Lie algebras.

\sekt{Dual pairs associated with nilpotent pairs\nopagebreak}%
Let $\ah,\ah'\subset\g$ be two subalgebras. Following R.\,Howe, we say 
that $\ah$ and $\ah'$ form a {\it dual pair}, if $\ah'=\z_\g(\ah)$
and vice versa.  A {\it reductive dual pair\/} is 
a dual pair $(\ah,\ah')$ such that each of $\ah,\ah'$ is reductive.
It is clear how to define a dual pair of groups.
In the group setting the problem is however more subtle, because of
connectedness questions.
A classification of reductive dual pairs in reductive Lie algebras
was obtained by H.\,Rubenthaler, see \cite{rub}. In the spirit of Dynkin,
he introduced the notion of an `$S$-irreducible' dual pair and described all
such pairs in the simple Lie algebras. The general classification is
then reduced to that for $S$-irreducible pairs. 
\begin{rem}{Definition} A dual pair $(\ah,\ah')$ is called 
$S$-{\it irreducible}, if
$\ah+\ah'$ is an $S$-subalgebra in the sense of Dynkin, i.e., it is not 
contained in a proper regular\footnote{A subalgebra of $\g$ is called 
{\it regular\/}
whenever its normalizer contains a Cartan subalgebra.} subalgebra of $\g$.
\end{rem}%
Let $\eb\in\fe$ be a nilpotent pair and $\hb$ 
a semisimple pair satisfying Eq.\,\re{comrel}. Then the quadruple $(\eb,\hb)$
is said to be {\it quasi-commutative}.
By definition, put $\ka_i=\z_\g(e_i,h_i)$, $i=1,2$. Our aim is to 
demonstrate a sufficient condition for $(\ka_1,\ka_2)$ to be a dual pair.
Note that $e_2,h_2\in\ka_1$ and $e_1,h_1\in\ka_2$. 
Consider the bi-grading of $\g$ determined by $\hb$:
\quad $\g=\bigoplus_{i,j}\g_{i,j}$, where $(i,j)$ runs over a finite subset of
$\Bbbk\times\Bbbk$ including (0,0), (1,0), and (0,1). 
The restriction of this bi-grading to either $\ka_1$ or $\ka_2$ 
gives ordinary gradings
$\ka_1=\oplus_{j}(\ka_1)_j$ and $\ka_2=\oplus_{i}(\ka_2)_i$, where
$(\ka_1)_j\subset\g_{0,j}$ and $(\ka_2)_i\subset\g_{i,0}$.
\begin{s}{Proposition}  \label{dual-int}
Let $(e_1,e_2,h_1,h_2)$ be a quasi-commutative quadruple.
Suppose  $\dim\z_\g(e_1,h_1,e_2)=\dim\z_\g(e_1,h_1,h_2)$.
Then \par
{\rm (i)} the grading of
$\ka_1$ 
is actually a $\Bbb Z$-grading, i.e., the
eigenvalues of $\ad h_2$ 
on $\ka_1$ 
are integral. Furthermore, the centralizer $\z_{\ka_1}(e_2)=\z_\g(e_1,h_1,e_2)$
is nonnegatively graded;  \par
{\rm (ii)} $(\ad e_2)_j : (\ka_1)_{j}\rar (\ka_1)_{j+1}$ is onto 
for $j\ge 0$.\\
{\ququ [Of course, this has the symmetric analogue, where indices 1 and 2
are interchanged.]}
\end{s}\begin{proof*}  (i)
The space $\z_\g(e_1, h_1,h_2)=\z_{\ka_1}(h_2)=(\ka_1)_0$ possesses the
$e_2$-filtration and $\mbox{lim\,}_{e_2}\z_{\ka_1}(h_2)\subset
\z_{\ka_1}(e_2)=\z_\g(e_1,h_1,e_2)$. It follows from the definition of 
$e_2$-limit that $\mbox{lim\,}_{e_2}\z_{\ka_1}(h_2)\subset
\oplus_{j\in {\Bbb Z}_{\ge 0}}(\ka_1)_j$. Furthermore, 
$\dim(\mbox{lim\,}_{e_2}\z_{\ka_1}(h_2))=\dim\z_{\ka_1}(h_2)$.
Under our assumption, this means that
$\mbox{lim\,}_{e_2}\z_\g(e_1,h_1,h_2)=\z_\g(e_1,h_1,e_2)$ and the eigenvalues
of $\ad h_2$ on $\z_\g(e_1,h_1,e_2)$ are nonnegative integers.
Assume that $(\ka_1)_j\ne 0$ for some $j\in\Bbbk\setminus\Bbb Z$. Since
$(\ka_1)_j$ is killed by some power of $\ad e_2$, we have
$j+c$ is the eigenvalue of $\ad h_2$ on $\z_\g(e_1,h_1,e_2)$
for some $c\in{\Bbb Z}_{\ge 0}$, which is impossible.
Thus, all the eigenvalues of $\ad h_2$ must be integral. \par
(ii) Set $(\ka_1)_{\ge j}=\oplus_{i\ge j}(\ka_1)_i$ and
consider the linear map 
$(\ad e_2)_{\ge 0} : (\ka_1)_{\ge 0}\rar (\ka_1)_{\ge 1}$.
By part (i), we have
$\mbox{Ker\,}(\ad e_2)_{\ge 0}=\z_{\ka_1}(e_2)$. 
That is, dimension of the kernel is
$\dim\z_\g(e_1,h_1,h_2)=\dim\,(\ka_1)_0$. Thus,
$(\ad e_2)_{\ge 0}$ must be onto.   
\end{proof*}%
\begin{s}{Theorem} \label{dual-main}
Suppose a quasi-commutative quadruple $(e_1,e_2,h_1,h_2)$ 
satisfies the conditions \par
1. $[\z_\g(e_1,h_1,h_2),\z_\g(e_2,h_1,h_2)]=0$,
\par
2. $\dim\z_\g(e_1,h_1,e_2)=\dim\z_\g(e_1,h_1,h_2)$, \par
3. $\dim\z_\g(e_2,h_2,e_1)=\dim\z_\g(e_2,h_2,h_1)$. \\
Then 
$(\ka_1,\ka_2)$ is a dual pair in $\g$. 
\end{s}\begin{proof}
Since $e_1,h_1\in\ka_2$ and $e_2,h_2\in\ka_1$, we have
$\ka_1\supset\z_\g(\ka_2)$ and $\ka_2\supset\z_\g(\ka_1)$. That is, the 
property of being a dual pair is
equivalent to that $[\ka_1,\ka_2]=0$. \\
We first prove
that $(\ka_2)_{\ge 0}$ commutes with $(\ka_1)_{\ge 0}$. 
Condition 1 says that $(\ka_1)_0$ commutes with $(\ka_2)_0$. 
Therefore the subalgebras generated by $\{(\ka_1)_0$, $e_2\}$ and 
$\{(\ka_2)_0, e_1\}$  commute.
By \re{dual-int}(ii), the subalgebra generated by $(\ka_1)_0$
and $e_2$ is $(\ka_1)_{\ge 0}$. Under condition 3, the same applies to 
$\ka_2$ in place of $\ka_1$.
That is, the subalgebra generated by $(\ka_2)_0$
and $e_1$ is $(\ka_2)_{\ge 0}$.
\\[.5ex]
Consider the
set ${\frak M}=\{[x,y]\mid x\in\ka_1, y\in\ka_2\}$. It is immediate that
${\frak M}$ is $\ad e_i$- and $\ad h_i$-stable $(i=1,2)$. Assume that ${\frak M}\ne\{0\}$,
that is, $[x,y]\ne 0$ for some $x\in (\ka_1)_j$ and $y\in (\ka_2)_i$.
By successively applying $\ad e_1$ and $\ad e_2$, we eventually obtain a
nonzero commutator $[x',y']$ with $x'\in (\ka_1)_{j'}$ and
$y'\in (\ka_2)_{i'}$ such that $x'\in\z_{\ka_1}(e_2)$ and
$y'\in\z_{\ka_2}(e_1)$.
It then follows from  \re{dual-int}(i) that $i'\ge 0$ and $j'\ge 0$. 
Thus, $x'\in (\ka_1)_{\ge 0}$, $y'\in (\ka_2)_{\ge 0}$ and one 
must have $[x',y']=0$.
This contradiction proves that ${\frak M}=\{0\}$.
\end{proof}%
Given $\eb$, it may {\sl a priori} happen that there are several non-equivalent
choices of $\hb$ such that $\hb$ satisfies Eq.\,\re{comrel} and the hypotheses
in \re{dual-main}. Fortunately, this question does not arise for (almost)
\pn pairs. We may even give a more precise statement in these cases.
Set $K_i:=Z_G(e_i,h_i)$, $i=1,2$. 
These groups  are not necessarily connected, but Lie$K_i=\ka_i$. 
%
\begin{s}{Theorem}  \label{dual-apn}
Suppose $\eb$ is either a \pn pair or an almost \pn pair and 
$\hb$ is an associated semisimple pair. Then \par
1. $(\ka_1,\ka_2)$ is a dual pair. The centre of $\ka_i$ ($i=1,2$) is 
trivial;   \par
2. This dual pair is reductive if and only if the pair $\eb$ is
rectangular;\par
3. $K_2=Z_G(K_1^o)$ and $K_1=Z_G(K_2^o)$; \par
4. If $\eb$ is a rectangular \pn pair, then $(\ka_1,\ka_2)$ is $S$-irreducible.
\end{s}\begin{proof*}
1. Since $\z_\g(\hb)$ is Abelian in both cases, hypothesis 1 in \re{dual-main}
is satisfied. By Theorem \ref{alm3}(3), the other hypotheses are satisfied, 
too.
The centre of $\ka_i$ is equal to $\ka_1\cap\ka_2=\z_\g(\eb)\cap\z_\g(\hb)=
\{0\}$.
\\[.5ex]
2. Clearly, $\ka_1$ is reductive if and only if $\ka_2$ is reductive.
If $\ka_2$ is reductive, then it contains a suitable
$\tri$-triple together with $e_1$. The opposite implication follows from
\ref{triple}(3).
\\[.5ex]
3. By symmetry, it suffices to prove the first equality.
Since $e_2,h_2\in\ka_1$, we have 
\[ 
K_2=Z_G(e_2,h_2)\supset Z_G(\ka_1)=Z_G(K_1^o)\ .
\]
In the proof of the opposite inclusion we use the relation 
$Z_G(\ka_1)\supset K_2^o$ proved in the first part.
Let $s\in K_2$ be an arbitrary element. One has to prove that $s{\cdot}x=x$ for
all $x\in\ka_1$. 
By \re{alm3}(4), $(\ka_1)_0$ is just the centre of
$\el_2$. Because $K_2$ lies in the {\it connected\/} group $L_2:=Z_G(h_2)$,
it commutes with $(\ka_1)_0$. By the very definition, $K_2$ commutes with
$e_2$. Thus, it commutes with $(\ka_1)_{\ge 0}$.
It then follows from \re{dual-int}(i) that $s{\cdot}x=x$ for $x\in\z_{\ka_1}(e_2)$. 
Consider
${\cal Y}:=\{y\in\ka_1\mid s{\cdot}y\ne y \}$. Suppose 
${\cal Y}\neq\varnothing$. Choose an element $y_0\in \cal Y$ which is killed
by the {\it least\/} possible power, say $p$, of $\ad e_2$. That is,
$(\ad e_2)^py_0\ne 0$ and $(\ad e_2)^{p+1}y_0=0$. Since 
${\cal Y}\cap\z_{\ka_1}(e_2)=\varnothing$, we have $p\ge 1$. Then
$(\ad e_2)^py_0\in\z_{\ka_1}(e_2)\subset (\ka_1)_{\ge 0}$ and 
hence $(\ad e_2)^py_0=
s{\cdot}((\ad e_2)^py_0)=(\ad e_2)^p(s{\cdot}y_0)$. In other words,
$(\ad e_2)^p(s{\cdot}y_0-y_0)=0$.
It follows that
$y_1:=s{\cdot}y_0-y_0\not\in{\cal Y}$ and $s{\cdot}y_1=y_1$. Therefore
$s^n{\cdot}y_0=y_0+ny_1$ for all $n\in \Bbb N$. However, we have
$s^n\in K_2^o\subset Z_G(\ka_1)$ for some $n>0$ and therefore 
$y_1$ must be zero.
This contradiction proves that ${\cal Y}=\varnothing$.
\\[.5ex]
4. It follows from parts 1 and 2 that $\ka_1+\ka_2$ is semisimple. Assume 
that $\ka_1+\ka_2\subset\g^{(1)}$, where $\g^{(1)}$ is a proper regular 
subalgebra of $\g$. Then there exists a maximal semisimple subalgebra 
$\g^{(2)}\subset\g^{(1)}$ such that $\ka_1+\ka_2\subset\g^{(2)}$. This
$\g^{(2)}$ is a regular subalgebra of $\g$, too. According to the description
of {\it maximal\/} regular semisimple subalgebras of $\g$, $\g^{(2)}$ is
contained in the fixed-point subalgebra of some element $s\in G$ of prime
order ($s\ne 1$). Then $s\in Z_G(\ka_1+\ka_2)$. However,
$Z_G(\ka_1+\ka_2)=Z_G(\eb)\cap Z_G(\hb)=\{1\}$, since $Z_G(\eb)$ is connected 
and unipotent \cite[3.6]{vitya}.
\end{proof*}
\\[.5ex]
{\bf Remark.} Arguing as in part 4 and using Prop.~\ref{con-apn}(1), 
one proves that if $\eb$ is either
a \pn pair or an almost \pn pair of $\Bbb Z$-type, then $\ka_1+\ka_2$ is
not contained in a proper {\it reductive\/} regular subalgebra of $\g$. 
However, $\ka_1+\ka_2$  
may lie in a proper parabolic subalgebra for a non-rectangular \pn pair $\eb$,
see example 1 in \re{primery}.
\begin{s}{Corollary}  \label{dualgroup}
If $K_1$ and $K_2$ are connected, then $(K_1,K_2)$ is a dual pair of groups
in $G$. \qu
\end{s}%
One may observe that we did not use properties of (almost) \pn pairs in
full strength in the above proofs. This suggests the notion of an (almost) \pn 
pair could be weakened so that the conclusion of Theorem
\ref{dual-main} still holds. 
A possible generalization in the rectangular case is discussed in the next
section.
\\[1ex]
{\bf Remark.} In the rectangular case, $K_i$ is a maximal reductive subgroup 
of $Z_G(e_i)$ and therefore 
$K_i/K_i^o\simeq Z_G(e_i)/Z_G(e_i)^o$. This group is 
known for all nilpotent orbits. The description is due to Springer and
Steinberg \cite{SS}
for the classical Lie algebras and due to A.\,Alekseevskii
\cite{andrei} 
for the exceptional ones. 

\begin{rem}{Examples}   \label{primery}
We give illustrations to theorem \re{dual-apn}.\\
{\sf 1}. The simplest non-rectangular \pn pair occurs in $\g={\frak sl}_3$. 
Let \\
$e_1{=}\left(\begin{array}{ccc} 0 & 0 & 1 \\
0 & 0 & 0 \\ 0 & 0 & 0 \end{array}\right)$ and
$e_2{=}\left(\begin{array}{ccc} 0 & 0 & 0 \\
0 & 0 & 1 \\ 0 & 0 & 0 \end{array}\right)$. Then $\z_\g(e_1,e_2)=\langle
e_1,e_2\rangle$, $h_1=\mbox{diag\,}(2/3,-1/3,-1/3)$ and 
$h_2=\mbox{diag\,}(-1/3,2/3,-1/3)$.
Whence $\ka_1=\langle e_2,h_2\rangle$ and $\ka_2=\langle e_1,h_1\rangle$.
It is clearly visible in this case that, for instance,
$\tih_1=\mbox{diag\,}(1,0,-1)$ and $\tih_1\ne 2h_1$. 
\\[.6ex]
{\sf 2}. In \re{seriya}, a series of almost \pn pairs in ${\frak sp}_{4n}$ is 
described.
In that case $\ka_1=\langle e_2,h_2\rangle$, while 
$\dim\ka_2=2n^2-n+1$. The Levi decomposition of $\ka_2$
is as follows: $\ka_2^{red}\simeq {\frak so}_{2n-1}\oplus\Bbbk$;
$\ka_2^{nil}$ is Abelian and affords the simplest representation of
${\frak so}_{2n-1}$.
\\[.6ex]
{\sf 3}. The rectangular \pn pairs in simple Lie algebras 
were classified  in \cite{wir}. For instance,
there are 4 such pairs in $\GR{E}{7}$ and
1 pair in either of $\GR{F}{4}$, $\GR{E}{6}$, $\GR{E}{8}$.
The corresponding $S$-irreducible
 reductive dual pairs are: \par
$(\GR{G}{2},\GR{A}{1})$ in $\GR{F}{4}$;
$(\GR{G}{2},\GR{A}{2})$ in $\GR{E}{6}$;
$(\GR{G}{2},\GR{F}{4})$ in $\GR{E}{8}$; \par
$(\GR{G}{2},\GR{A}{1})$, $(\GR{G}{2},\GR{C}{3})$,
$(\GR{F}{4},\GR{A}{1})$, $(\GR{A}{1},\GR{A}{1})$
 in $\GR{E}{7}$; \\[.5ex]
It turns out {\sl a posteriori\/} that for all nilpotent orbits 
occurring in this
situation the groups $Z_G(e_i)$ are connected (in the adjoint group!).
Hence the above two lines represent also the dual pairs of connected groups in
the respective adjoint group $G$.
\end{rem}%
{\bf Remark.}
One cannot hope to detect new instances of {\it reductive\/} dual pairs, 
since their classification was obtained by Rubenthaler. 
However, the `tableau r\'ecapitulatif' in \cite[p.\,70]{rub} contains several
inaccuracies. Below we use Rubenthaler's notation.
Each time an orthogonal Lie algebra $o(m)$ occurs as factor, one
has either to require that $m\ne 2$,  or to replace the given dual pair by a
correct one. This refers to the following possibilities in that table: \par
$\GR{B}{n}$ : $2(n-kp)+1-p=2$;\quad
$\GR{C}{n}$\,2) : $p=2$; \quad
$\GR{D}{n}$\,1) : $2n-2kp-p=2$. \\
For instance, if $p=2$ for $\GR{C}{n}$, then the dual pair must be
$({\frak gl}(k+1), o(2))$, not $({\frak sp}(k+1), o(2))$. However, unlike
the case $p\ne 2$, this dual pair is not $S$-irreducible. \\
It is also interesting to observe that Rubenthaler's ``diagrammes en dualit\'e"
correspond exactly to the dual pairs arising from the rectangular \pn pairs.


\sekt{Semi-principal pairs\nopagebreak}%
We shall say that a subalgebra $\ah\subset
\g$ is {\it reflexive} whenever $\z_\g(\z_\g(\ah))=\ah$. This is
tantamount to saying that $(\ah,\z_\g(\ah))$ is a dual pair. Obviously,
$\z_\g(\ah)$ is reflexive for {\it any\/} algebra $\ah\subset\g$. 
In particular, the centralizer of any $\tri$-triple is reflexive.
It is therefore interesting to to find out those $\tri$-triples whose double 
centralizer has some natural description, e.g.,
is again the centralizer of an $\tri$-triple. For instance,
in the dual pair associated to a rectangular \pn pair, both
algebras $\ka_1$ and $\ka_2$ are the centralizers of $\tri$-triples. 
Moreover, $\ka_1$ can be described as the centralizer in $\g$ of a
principal $\tri$-triple in $\ka_2$. This fact and the criterion given in
\re{th-1A} provide some motivation for the following definition. 
Recall that a nilpotent element $e$ 
in a reductive Lie algebra $\el$ 
is called {\it distinguished} whenever any semisimple element of
$\z_\el(e)$ lies in the centre of $\el$.
\begin{rem}{Definition}  \label{spnp}
A pair of nilpotent elements $(e_1,e_2)\in\g\times\g$ is called 
{\it semi-principal rectangular} (= \spr pair), if the following holds: \par
(i) there exist commuting $\tri$-triples 
$\{e_1,\tih_1,f_1\}$ and $\{e_2,\tih_2,f_2\}$ (rectangularity);  \par
(ii)  $e_1$ is distinguished in $\z_\g(\tih_2)=:\el_2$; \par
(iii)  $e_2$ is even in $\z_\g(e_1,\tih_1)=:\ka_1$.
\end{rem}%
Define the subalgebras $\ka_i$, $\el_i$, $\ce_i$ and the subgroups $K_i$
($i=1,2$) as above, with $\tih_i$ in place of $h_i$.
The meaning of condition (ii) is that $e_1$ should be a distinguished
element in $\ka_2$ which remains distinguished as element of $\el_2$.
Note that $(e_2,e_1)$ need not be
an \spr pair and $e_2$ need not be even in $\g$.
But if $e_2$ is even in $\g$, it is also even in $\ka_1$. It follows from
\re{alm3} that each rectangular (almost) \pn pair is an \spr pair. 
\begin{s}{Theorem}  \label{dualrect}
Let $(e_1,e_2)$ be an \spr pair. Then 
\par
{\rm (i)} $\z_\g(e_1,\tih_1,\tih_2)=\ce_2$;
\par
{\rm (ii)} $\ce_2$ is a Cartan subalgebra and  $e_2$ is a regular nilpotent
element in $\ka_1$;
\par
{\rm (iii)} $(\ka_1,\ka_2)$ is a reductive dual
pair in $\g$.
\end{s}%
\begin{proof*}
(i) The argument is close to that in \ref{alm3}(4).
By definition, $\z_\g(e_1,\tih_2)=\z_{\el_2}(e_1)=\ce_2\oplus\n$,
where $\n\subset [\el_2,\el_2]$ consists of nilpotent elements.
As $e_1,\tih_1\in\el_2$, we have $\z_\g(e_1,\tih_1,\tih_2)\supset\ce_2$.
Thus, $\ce_2\subset
\z_\g(e_1,\tih_1,\tih_2)\subset\z_{\g}(e_1,\tih_2)=\ce_2\oplus\n$.
In the rectangular case, $\ka_1=\z_\g(e_1,\tih_1,f_1)$ is reductive.
Hence
$\z_\g(e_1,\tih_1,\tih_2)=\z_{\ka_1}(\tih_2)$ is reductive, too. This clearly
forces that $\z_{\ka_1}(\tih_2)=\ce_2$.
\par
(ii) Since $\tih_2$ is semisimple, the previous equality means $\ce_2$ 
is a Cartan subalgebra
in $\ka_1$. Because $e_2$ is assumed to be even in $\ka_1$,
the $\tri$-triple $\{e_2,\tih_2,f_2\}$ is principal in $\ka_1$. \par
(iii) As in the proof of
\re{dual-main}, it is enough to prove that $[\ka_1,\ka_2]=0$.
It follows from (ii) that
$\ka_1$ is generated by $\ce_2, e_2$, and $f_2$ as Lie algebra. Since
$\ka_2\subset\el_2$ and $\ce_2$ is the centre of $\el_2$, we see that $\ka_2$ 
commutes with the stuff just described. 
\end{proof*}
\begin{rem}{Examples} \label{exc-primery}
A method of searching \spr pairs is as follows.
Let $\{e_2,\tih_2,f_2\}$ be an $\tri$-triple. First, one has to explicitly
determine $\ka_2$, $\el_2$, and the embedding $\ka_2\hookrightarrow\el_2$. 
A universal way is to exploit the weighted Dynkin diagram of $e_2$, but
in case of classical Lie algebras the formulas in terms of partitions are
also available. The next step is to find a distinguished element $e_1\in\ka_2$
which remains distinguished in $\el_2$. In case $e_2$ being even in $\g$, 
this is 
enough. Otherwise, one need to check that $e_2$ is even in $\ka_1$. 
The first candidate for $e_1$ is a \rn element in $\ka_2$. 
However, it can happen that \rn elements in $\ka_2$ fail to be 
distinguished in $\el_2$, while elements of a smaller orbit in $\ka_2$
satisfy our requirements. Furthermore, it can happen that there are several
such orbits in $\ka_2$. This means we obtain \spr pairs $(e',e_2)$ and 
$(e'',e_2)$ such that $e'$ and $e''$ lie in different $G$-orbits in $\g$,
see e.g. example 3 below.  Nevertheless,
it follows from \re{dualrect}(iii) that reductive parts of the centralizers of
$e'$ and $e''$ will coincide---they are just equal to $\z_\g(\ka_2)$. 
\\[.5ex]
We refer to \cite[ch.\,4\,\&\,8]{CoMc} for standard
facts on weighted Dynkin diagrams and labelling of nilpotent orbits.
\\[.5ex] \indent
{\sf 1}. Let $\co_2$ be the nilpotent orbit in $\g=\GR{E}{7}$, labelled by
$2\GR{A}{2}$. The weighted Dynkin diagram of $\co_2$ is
$\biggl($~
\begin{tabular}{@{}c@{}}
0--2--0--\lower3.45ex\vbox{\hbox{0\rule{0ex}{2.5ex}}
\hbox{\hspace{0.4ex}\rule{.1ex}{1ex}\rule{0ex}{1.4ex}}\hbox{0\strut}}--0--0
\end{tabular}
~$\biggr)$. 
Therefore $\el_2\simeq {\frak so}_{10}\oplus {\frak sl}_{2}\oplus\Bbbk$
and one finds in \cite{El75} that 
$\ka_2\simeq\GR{G}{2}\oplus{\frak sl}_{2}$. The embedding
$\ka_2\hookrightarrow [\el_2,\el_2]$ is  as follows:
$[\el_2,\el_2]$ has tautological 12-dimensional module
$\V_{12}=\V_{10}\oplus \V_2$. Then $\V_{10}\mid_{\ka_2}=
(\mbox{7-dim repr. }\GR{G}{2})\oplus \ad{\frak sl}_2$ and $\V_2\mid_{\ka_2}=
(\mbox{2-dim repr. }{\frak sl}_2)$.
Let $e_1$ be a \rn element in $\ka_2$. The above description of embedding shows
that $e_1$ is distinguished 
as element of $\el_2$. More precisely, $e_1=e'+e''$, where 
$e'\in{\frak so}_{10}$ corresponds to the partition (7,3) and 
$e''\in {\frak sl}_2$ is regular. (The distinguished nilpotent orbits in
${\frak so}_N$ correspond bijectively
to the partitions of $N$ into distinct odd parts.)  Since $e_2$ is even in
$\g$, it is also even in $\ka_1$. Hence a dual pair comes up and it 
remains to realize what $\ka_1$ is. 
The orbit of $e'$ in ${\frak so}_{10}$ is subregular and
is labelled by $\GR{D}{5}(a_1)$. Therefore the label
of $\co_1=G{\cdot}e_1$ is $\GR{D}{5}(a_1)+\GR{A}{1}$. Now, one finds in
the list of weighted Dynkin diagrams for $\GR{E}{7}$ that 
the diagram corresponding to $\co_1$ is 
$\biggl($~
\begin{tabular}{@{}c@{}}
0--0--2--\lower3.45ex\vbox{\hbox{0\rule{0ex}{2.5ex}}
\hbox{\hspace{0.4ex}\rule{.1ex}{1ex}\rule{0ex}{1.4ex}}\hbox{0\strut}}--0--2
\end{tabular}
~$\biggr)$. 
Hence $\el_1\simeq {\frak sl}_{4}\oplus {\frak sl}_{3}\oplus\Bbbk^2$
and, by \cite{El75}, $\ka_1\simeq{\frak sl}_{2}$.
Thus, the  dual pair 
is $(\GR{A}{1}, \GR{G}{2}\oplus \GR{A}{1})$. 
As in example~3 of \re{primery}, the groups $Z_G(e_i)$ ($i=1,2$) appear to 
be connected. The connected groups $K_i$ ($i=1,2$) form therefore a dual 
pair of groups.
\\[.6ex]
\indent
{\sf 2}. The members of \spr pairs are not necessarily even.
Let $\co_2$ be the nilpotent orbit in $\g=\GR{E}{7}$, labelled by
$2\GR{A}{1}$. Its weighted Dynkin diagram is
$\biggl($~
\begin{tabular}{@{}c@{}}
0--1--0--\lower3.45ex\vbox{\hbox{0\rule{0ex}{2.5ex}}
\hbox{\hspace{0.4ex}\rule{.1ex}{1ex}\rule{0ex}{1.4ex}}\hbox{0\strut}}--0--0
\end{tabular}
~$\biggr)$. 
Here $\el_2\simeq{\frak so}_{10}\oplus {\frak sl}_{2}\oplus\Bbbk$ and
$\ka_2\simeq{\frak so}_{9}\oplus{\frak sl}_{2}$ with the obvious embedding.
Therefore a \rn element $e_1\in\ka_2$ is also regular in  $\el_2$.
The label of $G{\cdot}e_1$ is $\GR{D}{5}+\GR{A}{1}$ and the 
weighted Dynkin diagram is
$\biggl($~
\begin{tabular}{@{}c@{}}
0--1--1--\lower3.45ex\vbox{\hbox{0\rule{0ex}{2.5ex}}
\hbox{\hspace{0.4ex}\rule{.1ex}{1ex}\rule{0ex}{1.4ex}}\hbox{1\strut}}--1--2
\end{tabular}
~$\biggr)$. Then one finds $\ka_1\simeq{\frak sl}_{2}$.
Hence $e_2$ is certainly even in $\ka_1$ and we obtain an \spr pair.
The corresponding reductive dual pair is $(\GR{B}{4}+\GR{A}{1},\GR{A}{1})$.
It is not $S$-irreducible, since it is contained in
the semisimple subalgebra of maximal rank $\GR{D}{6}+\GR{A}{1}\subset
\GR{E}{7}$.
\\[.6ex] \indent
{\sf 3}. Give a series of \spr pairs in a classical Lie algebra. Let 
$\g={\frak sp}_{2N}$ and let $G{\cdot}e_2$ be the orbit corresponding to the
partition $(\underbrace{m,\dots,m}_{2n},\underbrace{1,\dots,1}_{2l})$, where 
$m$ is odd $(m\ne 1)$ and $N=nm+l$. 
Since the parts have the same parity, $e_2$ is even.
Making use of the weighted Dynkin diagram, 
one finds that
$[\el_2,\el_2]=({\frak sl}_{2n})^{(m-1)/2}\oplus{\frak sp}_{2(n+l)}$ and
$\ka_2={\frak sp}_{2n}\oplus{\frak sp}_{2l}$. The embedding $\ka_2
\hookrightarrow\el_2$ is determined by the maps $\nu_1: \ka_2\to
{\frak sp}_{2(n+l)}$ and $\nu_2: \ka_2\to ({\frak sl}_{2n})^{(m-1)/2}$.
The first one is just the direct sum of matrices and the second one 
corresponds to the diagonal embedding 
${\frak sp}_{2n}\hookrightarrow {\frak sl}_{2n}\stackrel{\Delta}{\to}
({\frak sl}_{2n})^{(m-1)/2}$. Let us try to realize which elements
$e_1=e'+e''$ ($e'\in {\frak sp}_{2n}$, $e''\in {\frak sp}_{2l}$)
remain distinguished in $\el_2$. Since the only distinguished elements
in ${\frak sl}_N$ are the regular ones, $e'$ must be regular in 
${\frak sp}_{2n}$. This already guarantees us that $\nu_2(e_1)$ is regular
in $({\frak sl}_{2n})^{(m-1)/2}$. The orbits of distinguished elements in
${\frak sp}_{2N}$ corresponds bijectively to the partitions of $2N$ into
even unequal parts. Since $e'$ is already chosen, the partition of 
$\nu_1(e_1)$ has a part equal to $2n$. Thus, $e''$ must be a distinguished
element in ${\frak sp}_{2l}$ whose partition contains \un{no} parts equal 
to $2n$.
For instance, one may take $e''$ to be regular whenever $n\ne l$. In case
$n=l$, it is easy to see that a required partition exists if and only if
$n\not\in\{1,2\}$. Thus, \spr pairs come up if and only if
$(n,l)\not\in\{(1,1),(2,2)\}$, and the choice of $e_1$ is not 
unique in general. The partition of $e_1$ is either of
$(\underbrace{2n,\dots,2n}_{m},2l_1,\dots,2l_t)$, where $\sum_il_i=l$,
$l_i\ne l_j$, and $l_i\ne n$. For all such choices, $\ka_1$ is equal to
${\frak so}_m$ and we obtain the dual pair 
$({\frak sp}_{2n}\oplus{\frak sp}_{2l},{\frak so}_m)$. 
It follows from
\cite{rub} that these algebras form a dual pair even if $(n,l)=(1,1)$
or (2,2). But, for these `bad' values ${\frak so}_m$ has no 
interpretation as the centralizer in $\g$ of an
$\tri$-triple in $\ka_2$.
Observe also that one obtains a rectangular \pn pair, if $l=0$.
\end{rem}%
It is not hard to classify all even elements arising as the second
member of an \spr pair.
But instead of doing this, we introduce 
and classify more general elements in the next sections.

\sekt{Excellent elements and excellent sheets\nopagebreak}%
For an \spr pair $(e_1,e_2)$, Theorem \ref{dualrect}
says that $e_2$ is a regular nilpotent element in $\ka_1=\z_\g(\ka_2)=
\z_g(\z_g(e_2,\tih_2,f_2))$ and that $\ce_2=\z_g(\z_g(\tih_2))$ is the
centralizer of $\tih_2$ in $\ka_1$. That is, the functor of taking the double 
centralizer, applied to $\{e_2,\tih_2,f_2\}$, has nice properties. Our goal in 
this section is to axiomatize, investigate, and give applications of this
phenomenon.
\begin{rem}{Definition}  \label{excel} 
A nilpotent element $e$ is called {\it quasi-excellent}, if
$\dim\z_\g(\z_\g(\tih))=\rk \z_\g(\z_\g(e,\tih,f))$ for a (any)
$\tri$-triple $\{e,\tih,f\}$ containing $e$; $e$ is called 
{\it excellent}, if it is even and quasi-excellent. 
The same terminology applies to the $\tri$-triple itself.
\end{rem}%
We set $\ka=\z_\g(e,\tih,f)$, $K=Z_G(e,\tih,f)$, $\ka^\vee=\z_\g(\ka)$, and 
$\el=\z_\g(\tih)$. Then $\ce:=\z_g(\z_\g(\tih))$ is the centre of $\el$
and $(\ka,\ka^\vee)$ is a dual pair. We shall write $\z^2_\g(\cdot)$ in place
of $\z_\g(\z_\g(\cdot))$.
\\[.5ex]
{\bf Examples.} 1. If $e$ is distinguished in $\g$, then $\ka=0$ and
$\ka^\vee=\g$. But $\ce$ is a Cartan subalgebra if and only if $e$ is regular.
Hence an excellent distinguished element is regular. \\
2. If $(e_1,e_2)$ is an \spr pair, then $e_2$ is quasi-excellent. But the 
converse is not true. If $(n,l)=(1,1)$ or (2,2) in example 3 in 
\re{exc-primery}, then $e_2$ is excellent and cannot be included in an \spr
pair.

\begin{s}{Theorem}  \label{excel1}
Let $e$ be quasi-excellent. Then \par
1. $\ce$ is a Cartan subalgebra and $\{e,\tih,f\}$ is a principal
$\tri$-triple in $\ka^\vee$; \par
2. $K=Z_G(\ka^\vee)$.
\end{s}\begin{proof}
1. Since $\tih\in\{e,\tih,f\}$, taking the double centralizer gives
$\ce\subset\ka^\vee$. Whence $\ce$ is Cartan in $\ka^\vee$. Next,
$\z_{\ka^\vee}(\tih)=\z_\g(\tih)\cap\ka^\vee=\z_\g(\ce)\cap\ka^\vee=
\z_{\ka^\vee}(\ce)=\ce$, which means $\tih$ is regular in $\ka^\vee$. 
The centralizer of $\{e,\tih,f\}$ in $\ka^\vee$ is equal to $\ka\cap\ka^\vee$,
the centre of $\ka^\vee$. That is,
$e$ is distinguished in $\ka^\vee$. Since
any distinguished element is even (see e.g. \cite[ch.\,8]{CoMc}),
the assertion follows. \\
2. As $\{e,\tih,f\}\subset\ka^\vee$, we obtain $K\supset Z_G(\ka^\vee)$. 
In view of part 1, $\ka^\vee$ is generated by $e,\ce$, and $f$ as Lie algebra.
By definition, $K$ centralizes $e$ and $f$; and $K$ centralizes
$\ce$, because $\ce$ is the center of $\el$ and $K$ is contained in the
{\it connected\/} group $L=Z_G(\tih)$. Hence $K\subset Z_G(\ka^\vee)$.
\end{proof}%
Recall from section~2 the notions of $e$-filtration and $e$-limit, which
apply to any nilpotent element and any linear subspace of $\g$.
\begin{s}{Lemma}  \label{excel3}
Let $e\in\g$ be an arbitrary nilpotent element. Then \par
{\rm (i)} $[\lim_e \z^2_\g(\tih),\lim_e \z_\g(\tih)]=0$; \par
{\rm (ii)} If $e$ is even, then $\lim_e \z^2_\g(\tih)\subset \z^2_\g(e)$ and
$\dim\z^2_\g(\tih)\le\dim\z^2_\g(e)$.
\end{s}\begin{proof*}
Note first that $\z^2_\g(x)$ is the centre of $\z_\g(x)$ for any $x\in\g$.\par
1. By definition, the linear space ${\lim}_eM$ is generated by all
elements of the form $(\ad e)^ix$ ($x\in M$) that lie in $\z_\g(e)$.
Let $x\in\z^2_\g(\tih)$ and $y\in\z_\g(\tih)$. If 
$0\ne (\ad e)^i x\in \z_\g(e)$ and
$0\ne (\ad e)^j y\in \z_\g(e)$, then $[(\ad e)^ix, (\ad e)^j y]=
\frac{i!\,j!}{(i+j)!}(\ad e)^{i+j}[x,y]=0$. \par
2. If $e$ is even, then $\dim\z_\g(\tih)=\dim\z_\g(e)$. 
Since $\dim(\lim_e\z_\g(\tih))=\dim\z_\g(\tih)$, we conclude that 
$\lim_e\z_\g(\tih)=\z_\g(e)$. Hence $\lim_e\z^2_\g(\tih)\subset 
\z^2_\g(e)$, by the first claim.
\end{proof*} 
\begin{s}{Theorem}  \label{excel2}
Let $e$ be excellent. Then \par
1. $\ka$ and $\ka^\vee$ are semisimple; \par
2. $\z^2_\g(e)$ is the centralizer of $e$ in $\ka^\vee$.
\end{s}\begin{proof}
1. Let $\z$ be the centre of $\ka^\vee$. Then $\z\subset\ce$ and therefore
$[\z,\el]=0$. In case $e$ is even, $\g$ is generated by $e,\el$, and $f$ as Lie
algebra. Hence $\z$ is in the centre of $\g$, i.e., $\z=0$. \\
2. Since $\z^2_\g(e)\subset\z_{\ka^\vee}(e)$ and $\dim\z_{\ka^\vee}(e)=
\dim\z_{\ka^\vee}(\tih)$, the assertion
follows from the previous Lemma.
\end{proof}%
{\bf Example.} The properties in \ref{excel3}(ii) and \ref{excel2} need not 
hold, if $e$ is not even. Let $e$ be a nilpotent element in ${\frak sl}_5$
whose weighted Dynkin diagram is (1--1--1--1). Since $\tih$ is regular 
semisimple, $\z^2_\g(\tih)$ is a Cartan subalgebra and hence $e$ is 
quasi-excellent. However in this case $\ka$ is a 1-dimensional toral subalgebra
and $\dim\z^2_\g(e)=2<\rk\ka^\vee=4$. 
\\[.6ex]
We may express beautiful properties possessed by the excellent elements
(or excellent $\tri$-triples) in the following form.
If $\{e,\tih,f\}$ is excellent, then: \par
$\dim\z^2_\g(e)=\dim\z^2_\g(\tih)$; \par
$\z^2_\g(e)$ is the centalizer of $e$ in $\z^2_\g(e,\tih,f)$; \par
$\z^2_\g(\tih)$ is the centalizer of $\tih$ in $\z^2_\g(e,\tih,f)$.
\begin{subs}{Sheets} \end{subs}
Now, we show that the excellent elements provide an excellent
framework for constructing sections of sheets.
A {\it sheet} in $\g$ is an irreducible component of the set of points
whose $G$-orbits have a fixed dimension. The unique open sheet consists of
the regular elements in $\g$. This sheet has been thoroughly studied in
\cite{bert}.
The general theory of sheets 
was started in \cite{bk}. We refer to that paper 
for basic results of the theory. Each sheet is locally closed and contains
a unique nilpotent $G$-orbit. However, a nilpotent orbit may lie in several
sheets. We shall only deal with {\it Dixmier sheets\/}, i.e., sheets containing
semisimple elements. These are described as follows. For 
$Z\subset\g$, we set 
$Z^{reg}=\{x\in Z\mid \dim G{\cdot}x\ge\dim G{\cdot}y \mbox{ for all }
y\in Z\}$. Let $\el\subset\g$
be a Levi subalgebra with centre $\ce$. Then $(\ov{G{\cdot}\ce})^{reg}$
is a Dixmier sheet and all Dixmier sheets are of this form. 
To any even nilpotent element, one naturally associates a Dixmier sheet. 
If $e$ is even and $\tih$ is a characteristic of $e$,
then applying the above construction to 
the centre of 
$\z_\g(\tih)$, one obtains a Dixmier sheet containing $e$. This sheet will
be denoted by $\cs_{\tih}(e)$. In this case, one has
 $\dim\cs_{\tih}(e)=\dim G{\cdot}e+\dim\ce=\dim G{\cdot}\tih+\dim\ce$.
 Let us say that $Y\subset\cs_{\tih}(e)$ is a {\it section\/}
if $Y$ is irreducible,
$G{\cdot}y\cap\cs_{\tih}(e)=\{y\}$ 
for all $y\in Y$, and $G{\cdot}Y=\cs_{\tih}(e)$.
\\
In addition to the notation in \re{excel}, let $K^\vee$ denote the 
{\it connected\/}  group with Lie algebra $\ka^\vee$.
\begin{s}{Theorem}  \label{section}
Suppose $e$ is excellent and $(e,\tih,f)$ is an $\tri$-triple. Then \par
1.  $\cs_{\tih}(e)$ is smooth; \par
2. $\cs_{\tih}(e)\cap (e+\z_\g(f))=e+\z_{\ka^\vee}(f)$; \par
3. $e+\z_{\ka^\vee}(f)$ is a section of $\cs_{\tih}(e)$; \par
4. $\cs_{\tih}(e)$ is the unique sheet containing $e$.
\end{s}%
\begin{proof*} 1. Since $[\g,e]\oplus\z_\g(f)=\g$, the affine space
$e+\z_\g(f)$  is transversal to the orbit
$G{\cdot}e$ at $e$. Consider the subspace 
$\ca:=e+\z_{\ka^\vee}(f)=(e+\z_\g(f))\cap\ka^\vee$. Since
$\{e,\tih,f\}$ is a principal $\tri$-triple in $\ka^\vee$
(see \ref{excel1}(1)), $\ca$ is a
section of the open sheet in $\ka^\vee$. This is a classical result of
Kostant \cite{bert}. Therefore almost all elements in $\ca$ are 
semisimple and $K^\vee$-conjugate to 
elements in $\ce$, the latter being both a Cartan subalgebra in $\ka^\vee$
and the centre of $\el$. It follows that
$\max_{x\in\ca}\dim G{\cdot}x=\dim G-\dim \el=\dim\,G{\cdot}\tih$ 
and $\ov{G{\cdot}\ca}=
\ov{G{\cdot}\ce}=\ov{\cs_{\tih}(e)}$.
Consider the 1-parameter group $\{\lb(t)\mid t\in\Bbbk^*\}\subset GL(\g)$, 
where $\lb(t)=\mbox{exp}(t(\ad \tih-2{\cdot}\mbox{Id}_\g))$. It is easily
seen that $\ca$ is $\lb(\Bbbk^*)$-stable and
$e\in\ov{\lb(\Bbbk^*)x}$ for all $x\in\ca$.
Whence $\dim G{\cdot}e\le
\dim G{\cdot}x$. Because $e$ is assumed to be even and hence 
$\dim G{\cdot}e=\dim G{\cdot}\tih$, all $G$-orbits
intersecting $\ca$ have the same dimension.
Thus $\ca\subset\cs_{\tih}(e)$. 
\\[1ex]
Our next argument relies on results of P.\,Katsylo \cite{pasha}. He studied 
the variety $\cs\cap (e'+\z_\g(f'))$ for an arbitrary sheet $\cs$
containing an arbitrary nilpotent element $e'$.
By \cite[0.1]{pasha}, we have \par 
$\bullet$ $\cb:=\cs\cap (e'+\z_{\g}(f'))$ is closed in 
$e'+\z_{\g}(f')$, \par
$\bullet$ the $G$-orbits in $\cs$ intersect $\cb$
transversally, \par
$\bullet$ $G{\cdot}\cb_i=\cs$ for any 
irreducible component $\cb_i$ of $\cb$. 
\\[0.5ex]
Applying this to $\cb^{\tih}:=\cs_{\tih}(e)\cap (e+\z_{\g}(f))$
and the irreducible components $\cb^{\tih}_i$, we see that
 $\dim\cb^{\tih}_i=\dim\ce$.
Since $\dim\ca=\dim\z_{\ka^\vee}(f)=\dim\ce$ 
and $\ca\subset\cb^{\tih}$, we have $\ca$ is an irreducible component 
of $\cb^{\tih}$ and $\cs_{\tih}(e)= G{\cdot}\ca$.
It follows from transversality condition that the natural map
$G\times\ca\rar\cs_{\tih}(e)$ is smooth and hence $\cs_{\tih}(e)$ is
smooth, too.
\par
2. By \cite[0.2]{pasha}, 
the connected group $K^o$ acts trivially on $\cb^{\tih}$ or, equivalently, 
$\cb^{\tih}$ is contained in $\z_\g(\ka)=\ka^\vee$. Therefore
\[
\ca=(e+\z_\g(f))\cap\ka^\vee\subset(e+\z_\g(f))\cap\cs_{\tih}(e)
=\cb^{\tih}\subset\ka^\vee \ .
\]
Whence $\ca=\cb^{\tih}$. \par
3. By \cite[0.3]{pasha}, two points
$x',x''\in\ca$ lie in the same $G$-orbit if and
only if these lie in the same $K/K^o$-orbit.
Thus,
$\ca$ is a section of $\cs_{\tih}(e)$ if and only if $K$ acts
trivially on $\ca$. Let $x'$ be a generic point in $\ca$. Then $x'$ is
a regular semisimple element in $\ka^\vee$ and hence $K^\vee{\cdot}x'$ contains
a point $y\in\ce$. We have $Z_G(y)=Z_G(\tih)\supset K$ and
 $x'=s{\cdot}y$  for some $s\in K^\vee$. Then 
$Z_G(x')\supset s{}Ks^{-1}$. By \ref{excel1}(2), the subgroups $K$ and
$K^\vee$ commute. Hence $K\subset Z_G(x')$ and we are done.
\par
4. Let $\cs$ be an arbitrary sheet containing $e$. Arguing as in the proof 
of part 2, we obtain $(e+\z_\g(f))\cap\cs\subset\ka^\vee$. Therefore 
$(e+\z_\g(f))\cap\cs\subset e+\z_{\ka^\vee}(f)\subset\cs_{\tih}(e)$.
Since $\cs=G{\cdot}\Bigl( (e+\z_\g(f))\cap\cs\Bigr)$ by Katsylo's result, 
we must have
$\cs=\cs_{\tih}(e)$.
\end{proof*}
\begin{s}{Corollary}  \label{sectrect}
The assertions of \re{section} are valid for both members of 
the rectangular \pn pairs.
\end{s}%
\begin{proof} By theorems \ref{th-1A} and \ref{dualrect}, 
each member of a rectangular \pn pair
is excellent.
\end{proof}%
In view of \ref{section}(4), the sheet containing an excellent
element is said to be {\it excellent}, too. 
\\[.7ex]
One may remember that each sheet in ${\frak sl}_N$ is smooth and has 
a section, and each nilpotent element belongs to a unique sheet. On the
other hand, it is shown by Ginzburg that \pn pairs in ${\frak sl}_N$
are essentially being classified by Young diagrams with $N$ boxes, see
\cite[5.6]{vitya}. In particular, each nilpotent element 
can be included in a \pn pair. 
(This nice property is no longer true for the other simple
Lie algebras.) It is then natural to suggest that something like
Theorem \ref{section} holds for arbitrary \pn pairs. 
\begin{s}{Conjecture} \label{gip}
Let $e$ be a member of a \pn pair. Then $e$ belongs to a unique
sheet, this sheet is smooth and has a section.
\end{s}%
In principle, uniqueness of the sheet can be checked in a case-by-case
fashion.
The explicit description of induced orbits in the simple Lie algebras is 
known, see \cite{Sp}, \cite{ke83}, \cite{El85}\footnote{The results of this 
paper were announced in \cite[pp.\,171--177]{Sp}.}. Therefore, given a nilpotent 
orbit, one can say whether it belongs to a unique sheet.
Unfortunately, $\tri$-framework
breaks completely down in the general situation and it is not clear how to
produce a section. Some hint is however provided by \re{main1}:
If $(e_1,e_2)$ is a \pn pair, then $e_2$ is a Richardson element in $\p_2$.
Hence $e_2$ lies in the Dixmier sheet $(\ov{G{\cdot}\ce_2})^{reg}$. It
follows from \ref{smes}(4) and \ref{alm3}(3) that
$\dim\z_\g(e_1,h_1,e_2)=\dim\ce_2$ and one may hope that $\z_\g(e_1,h_1,e_2)$
should enter somehow in the conjectured section. At least, $f_2
+\z_\g(e_1,h_1,e_2)$ is a section in the rectangular case, cf. \re{section}.
Unfortunately, we have $2h_1\ne \tih_1$ and $f_2\not\in\z_\g(e_1,h_1)$ 
in general!

\sekt{Classification and tables  \nopagebreak}%
Since the excellent orbits (or sheets) enjoy excellent properties,
it is worth to get the list of them. Our
classification is presented in two tables.
Give the necessary details concerning our computations.

\begin{subs}{The exceptional case} \end{subs}
In $\GR{G}{2}$, the only excellent orbit is the regular nilpotent one. 
For the non-regular excellent
orbits, $\ka$ has to be non-trivial and semisimple.
Looking through the tables in \cite{El75}, one finds that the 
number of such even orbits in
$\GR{F}{4}$, $\GR{E}{6}$, $\GR{E}{7}$, $\GR{E}{8}$ is equal to 
3,\,3,\,15,\,13 respectively. Having computed $\z_\g(\ka)$ in each case, 
one distinguishes the excellent orbits among them.
The actual number of non-regular excellent orbits is equal to
2,\,2,\,9,\,6 respectively. 

\begin{subs}{The classical case} \end{subs}
If $d_1,\dots,d_m$ are all nonzero {\it different\/} parts of a partition 
$\db$ such that $d_1>d_2>\dots>d_m$
and $d_i$ occurs with multiplicity $r_i$ ($i=1,\dots,m$),
then we write $\db=({d_1}^{r_1},\dots,{d_m}^{r_m})$. 
For a classical simple Lie algebra, let $G(\db)$ denote the orbit 
corresponding to $\db$. It is assumed that $\db$ satisfies the necessary 
constraints in the symplectic and orthogonal case. (If $\g={\frak so}_N$ and
$\db$ is ``very even", then $G(\db)$ can be either of the two $SO_N$-orbits.)
It is well known (and easy to prove) that $G(\db)$ is even if and only if 
the $d_i$'s have the same parity.
\\[.5ex]
Given an even orbit $G(\db)$, we describe 
the structure of $\el$, $\ka$, and $\z_\g(\ka)$ in terms of $\db$.
The formulas for $\el$ are easy and those for $\ka$ are found in 
\cite[6.1.3]{CoMc}. Then it is not hard to realize what $\z_\g(\ka)$ is.
Some accuracy is however needed while dealing with algebras
${\frak so}_r$, since these are not semisimple for $r=2$. Since $\ka$ and
$\z_\g(\ka)$ must be semisimple for the excellent elements 
(see \ref{excel2}(1)), we will assume that
$r\ne 2$ whenever $\ka$ contains a summand ${\frak so}_{r}$. With explicit 
formulas for $\el$ and $\z_\g(\ka)$, verification of the arithmetical condition
$\dim\ce=\rk\z_\g(\ka)$ becomes trivial. \\
It is important to stress that our formulas for $\el$ are only valid for
\un{even} orbits.
\\[.6ex]
{\sf 1.} $\g={\frak sl}_N$. Here 
\[
\el=({\frak sl}_{r_1})^{d_1-d_2}\oplus 
({\frak sl}_{r_1+r_2})^{d_2-d_3}\oplus\dots\oplus
({\frak sl}_{r_1+\dots+r_m})^{d_m}\oplus \Bbbk^{d_1-1} \ ,
\]
\[ 
\ka={\frak sl}_{r_1}\oplus{\frak sl}_{r_2}\oplus\dots
\oplus{\frak sl}_{r_m}\oplus\Bbbk^{m-1} \ ,
\]
\[  \quad \mbox{and}\quad
\z_\g(\ka)={\frak sl}_{d_1}\oplus{\frak sl}_{d_2}\oplus\dots
\oplus{\frak sl}_{d_m}\oplus\Bbbk^{m-1} \ .
\]
Thus, $\ka$ is semisimple if and only if $m=1$ and then $G(\db)$ is 
excellent. Actually, $e$ is a member of a \pn pair in this case.
\\[.6ex]
{\sf 2.} $\g={\frak sp}_{2N}$. 
Now we have to distinguish two possibilities. 
\par (a) $d_1,\dots,d_m$ are odd. Then $r_1,\dots,r_m$ must be even.
Here
\[
\el=({\frak sl}_{r_1})^{(d_1-d_2)/2}\oplus 
({\frak sl}_{r_1+r_2})^{(d_2-d_3)/2}\oplus\dots\oplus
({\frak sl}_{r_1+\dots+r_m})^{(d_m-1)/2}\oplus{\frak sp}_{r_1+\dots+r_m}
\oplus\Bbbk^{(d_1-1)/2} \ ,
\]
%
\[
\ka={\frak sp}_{r_1}\oplus{\frak sp}_{r_2}\oplus\dots
\oplus{\frak sp}_{r_m} \ ,
\]
\[ \quad \mbox{and}\quad
\z_\g(\ka)={\frak so}_{d_1}\oplus{\frak so}_{d_2}\oplus\dots
\oplus{\frak so}_{d_m} \ .
\]
The arithmetical condition reads
$\frac{d_1-1}{2}=\frac{d_1-1}{2}+\dots+\frac{d_m-1}{2}$. Whence the excellent
orbits corresponds to either $m=1$ or
$m=2$ and $d_2=1$. The first possibility gives us a member of a \pn pair.
\par 
(b) $d_1,\dots,d_m$ are even. Then 
\[
\el=({\frak sl}_{r_1})^{(d_1-d_2)/2}\oplus 
({\frak sl}_{r_1+r_2})^{(d_2-d_3)/2}\oplus\dots\oplus
({\frak sl}_{r_1+\dots+r_m})^{d_m/2}
\oplus\Bbbk^{d_1/2} \ ,
\]
\[ \quad \mbox{and}\quad
\ka={\frak so}_{r_1}\oplus{\frak so}_{r_2}\oplus\dots
\oplus{\frak so}_{r_m} \ .
\]
If none of the $r_i$'s is equal to 2, then
$\displaystyle
\z_\g(\ka)=(\bigoplus_{i:\ r_i\ne 1}{\frak sp}_{d_i})\oplus
{\frak sp}_{d}$, where $\displaystyle d=\sum_{j:\ r_j=1}d_j$. But the rank of
$\z_\g(\ka)$ does not depend on the number of $r_i$'s that are equal to 1
and the condition 
$\frac{d_1}{2}=\frac{d_1}{2}+\dots+\frac{d_m}{2}$ implies that $m=1$.
\\[.6ex]
{\sf 3.} $\g={\frak so}_{N}$. Here again two possibilities.
\par (a) $d_1,\dots,d_m$ are even. Then $r_1,\dots,r_m$ must be even. Here
\[
\el=({\frak sl}_{r_1})^{(d_1-d_2)/2}\oplus 
({\frak sl}_{r_1+r_2})^{(d_2-d_3)/2}\oplus\dots\oplus
({\frak sl}_{r_1+\dots+r_m})^{d_m/2}\oplus\Bbbk^{d_1/2} \ ,
\]
\[
\ka={\frak sp}_{r_1}\oplus{\frak sp}_{r_2}\oplus\dots
\oplus{\frak sp}_{r_m} \ ,
\]
\[  \quad \mbox{and}\quad
\z_\g(\ka)={\frak sp}_{d_1}\oplus{\frak sp}_{d_2}\oplus\dots
\oplus{\frak sp}_{d_m} \ .
\]
Then the equality $d_1/2=(d_1+\dots+d_m)/2$ leads to $m=1$.
\par (b) $d_1,\dots,d_m$ are odd. Then 
\[
\el=({\frak sl}_{r_1})^{(d_1-d_2)/2}\oplus 
({\frak sl}_{r_1+r_2})^{(d_2-d_3)/2}\oplus\dots\oplus
({\frak sl}_{r_1+\dots+r_m})^{(d_m-1)/2}\oplus{\frak so}_{r_1+\dots+r_m}
\oplus\Bbbk^{(d_1-1)/2} \ ,
\]
\[  \quad \mbox{and}\quad
\ka={\frak so}_{r_1}\oplus{\frak so}_{r_2}\oplus\dots
\oplus{\frak so}_{r_m} \ .
\]
If none of the $r_i$'s is equal to 2, then
$\displaystyle
\z_\g(\ka)=(\bigoplus_{i:\ r_i\ne 1}{\frak so}_{d_i})\oplus
{\frak so}_{d}$, where $\displaystyle d=\sum_{j:\ r_j=1}d_j$.
Observe that an anomaly occurs if $r_1+\dots+r_m=2$, i.e.,
$m=2$ and $r_1=r_2=1$. Then $\dim\ce=\frac{d_1-1}{2}+1$. This case leads to 
the ``excellent" partition $(d_1,1)$, which represents the regular nilpotent 
orbit
in ${\frak so}_{d_1+1}$. Otherwise, we have $\dim\ce=\frac{d_1-1}{2}$ and
$\rk\z_\g(\ka)\ge \sum_i \frac{d_i-1}{2}$. Then a quick analysis leads to the 
following ``excellent" partitions:
$m=1$ and $r_1\ne 2$; $m=2$, $d_2=1$, and $r_i\ne 2$ $(i=1,2)$.
\\[.6ex] 
Thus, a classification of excellent orbits is completed.

\begin{table}[htbp]
\caption{The non-regular excellent orbits in the exceptional case} 
\label{vserect}
\vskip.8ex\centerline
{
\begin{tabular}{cccccc}
$\g$ & diagram of $G{\cdot}e$ & label of $G{\cdot}e$ & 
$[\el,\el]$ & $\ka$ & $\dim{\cal A}$  \\ \hline
$\GR{F}{4}$ & 2--0$\Leftarrow$0--0 & $\GRt{A}{2}$ & $\GR{B}{3}$ & 
$\GR{G}{2}$ & 1\\
& 0--0$\Leftarrow$2--2 & $\GR{B}{3}$ & $\GRt{A}{2}$ & $\GR{A}{1}$ & 2\\ \hline
$\GR{E}{6}$ & \begin{tabular}{@{}c@{}}
2--0--\lower3.3ex\vbox{\hbox{0\rule{0ex}{2.4ex}}
\hbox{\hspace{0.4ex}\rule{.1ex}{1ex}\rule{0ex}{1.1ex}}\hbox{0\strut}}--0--2
\end{tabular} 
 & $2\GR{A}{2}$  & $\GR{D}{4}$ & $\GR{G}{2}$ & 2 \\ 
 & \begin{tabular}{@{}c@{}}
0--0--\lower3.3ex\vbox{\hbox{2\rule{0ex}{2.4ex}}
\hbox{\hspace{0.4ex}\rule{.1ex}{1ex}\rule{0ex}{1.1ex}}\hbox{2\strut}}--0--0
\end{tabular} 
 & $\GR{D}{4}$  & $2\GR{A}{2}$ & $\GR{A}{2}$ & 2 \\ \hline
$\GR{E}{7}$ & \begin{tabular}{@{}c@{}}
2--0--0--\lower3.3ex\vbox{\hbox{0\rule{0ex}{2.4ex}}
\hbox{\hspace{0.4ex}\rule{.1ex}{1ex}\rule{0ex}{1.1ex}}\hbox{0\strut}}--0--0
\end{tabular} 
 & $[3\GR{A}{1}]''$  & $\GR{E}{6}$ & $\GR{F}{4}$ & 1 \\ 
& \begin{tabular}{@{}c@{}}
0--2--0--\lower3.3ex\vbox{\hbox{2\rule{0ex}{2.4ex}}
\hbox{\hspace{0.4ex}\rule{.1ex}{1ex}\rule{0ex}{1.1ex}}\hbox{0\strut}}--2--2
\end{tabular} 
 & $\GR{E}{6}$  & $[3\GR{A}{1}]''$ & $\GR{A}{1}$ & 4 \\ \cline{2-6} 
& \begin{tabular}{@{}c@{}}
0--0--0--\lower3.3ex\vbox{\hbox{0\rule{0ex}{2.4ex}}
\hbox{\hspace{0.4ex}\rule{.1ex}{1ex}\rule{0ex}{1.1ex}}\hbox{2\strut}}--0--0
\end{tabular} 
 & $\GR{A}{2}+3\GR{A}{1}$  & $\GR{A}{6}$ & $\GR{G}{2}$ & 1 \\ 
& \begin{tabular}{@{}c@{}}
0--2--0--\lower3.3ex\vbox{\hbox{2\rule{0ex}{2.4ex}}
\hbox{\hspace{0.4ex}\rule{.1ex}{1ex}\rule{0ex}{1.1ex}}\hbox{0\strut}}--0--0
\end{tabular} 
 & $\GR{A}{6}$  & $\GR{A}{2}{+}3\GR{A}{1}$ & $\GR{A}{1}$ & 2 \\ \cline{2-6}
& \begin{tabular}{@{}c@{}}
0--0--0--\lower3.3ex\vbox{\hbox{0\rule{0ex}{2.4ex}}
\hbox{\hspace{0.4ex}\rule{.1ex}{1ex}\rule{0ex}{1.1ex}}\hbox{0\strut}}--2--2
\end{tabular} 
 & $\GR{D}{4}$  & $[\GR{A}{5}]''$ & $\GR{C}{3}$ & 2 \\ 
& \begin{tabular}{@{}c@{}}
2--2--0--\lower3.3ex\vbox{\hbox{0\rule{0ex}{2.4ex}}
\hbox{\hspace{0.4ex}\rule{.1ex}{1ex}\rule{0ex}{1.1ex}}\hbox{0\strut}}--0--2
\end{tabular} 
 & $[\GR{A}{5}]''$  & $\GR{D}{4}$ & $\GR{G}{2}$ & 3 \\ \cline{2-6}
& \begin{tabular}{@{}c@{}}
0--0--2--\lower3.4ex\vbox{\hbox{0\rule{0ex}{2.4ex}}
\hbox{\hspace{0.4ex}\rule{.1ex}{1ex}\rule{0ex}{1.2ex}}\hbox{0\strut}}--0--0
\end{tabular} 
 & $\GR{A}{3}{+}\GR{A}{2}{+}\GR{A}{1}$  & $\GR{A}{4}{+}\GR{A}{2}$ & 
$\GR{A}{1}$ & 1 \\ 
& \begin{tabular}{@{}c@{}}
0--0--0--\lower3.4ex\vbox{\hbox{2\rule{0ex}{2.4ex}}
\hbox{\hspace{0.4ex}\rule{.1ex}{1ex}\rule{0ex}{1.2ex}}\hbox{0\strut}}--0--0
\end{tabular} 
 & $\GR{A}{4}+\GR{A}{2}$  & $\GR{A}{3}{+}\GR{A}{2}{+}\GR{A}{1}$ & 
 $\GR{A}{1}$ & 1 \\ \cline{2-6}
 & \begin{tabular}{@{}c@{}}
0--2--0--\lower3.4ex\vbox{\hbox{0\rule{0ex}{2.4ex}}
\hbox{\hspace{0.4ex}\rule{.1ex}{1ex}\rule{0ex}{1.2ex}}\hbox{0\strut}}--0--0
\end{tabular} 
 & ${2\GR{A}{2}}$ & $\GR{D}{5}{+}\GR{A}{1}$ 
& $\GR{G}{2}{+}\GR{A}{1}$ & 1\\ \hline
$\GR{E}{8}$  & \begin{tabular}{@{}c@{}}
2--2--0--0--\lower3.3ex\vbox{\hbox{0\rule{0ex}{2.4ex}}
\hbox{\hspace{0.4ex}\rule{.1ex}{1ex}\rule{0ex}{1.1ex}}\hbox{0\strut}}--0--0
\end{tabular} 
 & $\GR{D}{4}$ & $\GR{E}{6}$ & $\GR{F}{4}$ & 2 \\
 & \begin{tabular}{@{}c@{}}
2--2--2--0--\lower3.3ex\vbox{\hbox{0\rule{0ex}{2.4ex}}
\hbox{\hspace{0.4ex}\rule{.1ex}{1ex}\rule{0ex}{1.1ex}}\hbox{0\strut}}--0--2
\end{tabular} 
 & $\GR{E}{6}$ & $\GR{D}{4}$ & $\GR{G}{2}$ & 4 \\ \cline{2-6}
 & \begin{tabular}{@{}c@{}}
0--0--0--0--\lower3.4ex\vbox{\hbox{0\rule{0ex}{2.4ex}}
\hbox{\hspace{0.4ex}\rule{.1ex}{1ex}\rule{0ex}{1.2ex}}\hbox{0\strut}}--0--2
\end{tabular} 
 & ${2\GR{A}{2}}$ & $\GR{D}{7}$
& $2\GR{G}{2}$ & 1\\ \cline{2-6}
 & \begin{tabular}{@{}c@{}}
0--0--0--0--\lower3.4ex\vbox{\hbox{0\rule{0ex}{2.4ex}}
\hbox{\hspace{0.4ex}\rule{.1ex}{1ex}\rule{0ex}{1.2ex}}\hbox{2\strut}}--0--0
\end{tabular} 
 & $\GR{D}{4}(a_1){+}\GR{A}{2}$ & $\GR{A}{7}$
& $\GR{A}{2}$ & 1\\ \cline{2-6}
 & \begin{tabular}{@{}c@{}}
0--0--2--0--\lower3.4ex\vbox{\hbox{0\rule{0ex}{2.4ex}}
\hbox{\hspace{0.4ex}\rule{.1ex}{1ex}\rule{0ex}{1.2ex}}\hbox{0\strut}}--0--0
\end{tabular} 
 & $\GR{A}{4}{+}\GR{A}{2}$ & $\GR{D}{5}{+}\GR{A}{2}$
& $2\GR{A}{1}$ & 1\\ \cline{2-6}
 & \begin{tabular}{@{}c@{}}
0--0--2--0--\lower3.4ex\vbox{\hbox{0\rule{0ex}{2.4ex}}
\hbox{\hspace{0.4ex}\rule{.1ex}{1ex}\rule{0ex}{1.2ex}}\hbox{0\strut}}--0--2
\end{tabular} 
 & $\GR{A}{6}$ & $\GR{D}{4}{+}\GR{A}{2}$
& $2\GR{A}{1}$ & 2\\ \hline
\end{tabular}
}
\end{table}%
\noindent
Here are some explanations to the tables.
Nilpotent orbits in the exceptional (resp. classical) Lie
algebras are represented by their weighted Dynkin diagrams (resp. partitions).
The rightmost column gives dimension of the section of the 
excellent sheet. Recall from \re{section} that also $\dim\ca=\rk\z_\g(\ka)$.
In table \ref{vserect}, the pairs of orbits corresponding to the rectangular
\pn pairs are placed in adjacent rows that are not separated
 by horizontal line. The "duality"  between the label of $G{\cdot}e$
and the Cartan type of $[\el,\el]$ visible in each such pair is a 
manifestation of the properties stated in \re{main1} or in \re{th-1A}.
In table 2, the label of an orbit has the same meaning as for exceptional
Lie algebras. It represents the (unique up to conjugation) minimal Levi
subalgebra meeting this orbit. An algorithm for finding the label through
the partition is found in \cite[sect.\,3]{sph}.
\begin{table}[htbp]
\caption{The classical case} \label{klassrect}
\vskip1.5ex\centerline
{
\begin{tabular}{cccccc}
$\g$ & partition  & label of $G{\cdot}e$ & 
$[\el,\el]$ & $\ka$ & $\dim{\cal A}$  \\ \hline
${\frak sl}_{nm}$ 
  & $(n,...,n)$ 
  & $m\GR{A}{n-1}$ 
  & $({\frak sl}_{m})^n$ 
  & ${\frak sl}_{m}$ 
  & $n-1$  \\ \hline
$\begin{array}{c} {\frak sp}_{2nm}  \\
          \mbox{{\footnotesize ($m\ne 2$)}}
       \end{array}$
   & $(2n,...,2n)$ 
   & $\begin{array}{c} \frac{m-1}{2}\GRt{A}{2n-1}{+}\GR{C}{n} \\
                       \mbox{\footnotesize  if $m$ is odd;} \\
                       \frac{m}{2}\GRt{A}{2n-1} \\
                       \mbox{\footnotesize  if $m$ is even;}
      \end{array}$ 
   & $({\frak sl}_{m})^n$ 
   & ${\frak so}_{m}$ 
   & $n$  \\ \cline{2-6}
$\begin{array}{c} {\frak sp}_{2(nm+l)} \\ 
          \mbox{{\footnotesize ($m$ is odd)}}
       \end{array}$
  & $(m^{2n},1^{2l})$
  & $n\GRt{A}{m-1}$ 
  & $({\frak sl}_{2n})^{\frac{m-1}{2}}{\oplus}{\frak sp}_{2(n+l)}$ 
  & ${\frak sp}_{2n}{\oplus}{\frak sp}_{2l}$ 
  & $\frac{m-1}{2}$  \\ \hline
$\begin{array}{c} {\frak so}_{nm} \\
           \mbox{{\footnotesize ($m,n$ are even)}}
       \end{array}$
  & $(m,...,m)$ 
  & $\frac{n}{2}\GR{A}{m-1}$ 
  & $({\frak sl}_{n})^{m/2}$ 
  & ${\frak sp}_{n}$ 
  & $\frac{m}{2}$  \\ \cline{2-6}
$\begin{array}{c} {\frak so}_{nm+l} \\ 
             \mbox{\footnotesize ($m$ is odd)} \\ 
             \mbox{\footnotesize $n\ne 2, l\ne 2$}
  \end{array}$
   & $(m^n,1^l)$
   & $\begin{array}{c} \frac{n-1}{2}\GR{A}{m-1}{+}\GR{D}{\frac{m+1}{2}} \\
                       \mbox{\footnotesize  if $l$ is odd;} \\
                       \frac{n-1}{2}\GR{A}{m-1}{+}\GR{B}{\frac{m-1}{2}} \\ 
                       \mbox{\footnotesize  if $l$ is even;}
      \end{array}$ 
   & $({\frak sl}_{n})^{\frac{m-1}{2}}{\oplus}{\frak so}_{n+l}$ 
   & ${\frak so}_{n}{\oplus}{\frak so}_{l}$ 
   & $\frac{m-1}{2}$  \\ \hline
\end{tabular}
}
\end{table}%
\begin{rem}{Remark}  \label{Rubenth}
In \cite{nappes}, Rubenthaler introduced the notion of an {\it admissible\/} 
sheet
and proved that each admissible sheet has a section, which is an affine space.
He also gave a classification of the admissible sheets.
It follows from comparing the two classifications that 
each excellent sheet is admissible. But the converse is not true and,
furthermore, the assertions of Theorem \ref{section} do not hold for the 
nilpotent orbit lying in an arbitrary
admissible sheet. For instance, the nilpotent orbit labelled
by $\GR{D}{4}(a_1)$ in $\g=\GR{E}{6}$ lies in an admissible sheet, while the
total number of sheets containing it is equal to 3, see \cite[table 1]{El85}.
It should also be noted that Rubenthaler writes nothing about smoothness of
admissible sheets and that our approach to the problem is less technical.
\end{rem}

\vno{3} \indent
{\footnotesize 
\parbox{195pt}{%
{\it Math. Department \\
MIREA \\
prosp. Vernadskogo, 78 \\
Moscow 117454 \quad Russia} \\ panyush@dpa.msk.ru }
\hspace{1cm}
\parbox{190pt}{%
Current address (until April 28, 1999): \\ 
{\it Universit\'e de Poitiers \\
Math\'ematiques \quad S.P.2M.I. \\
Boulevard 3 -- T\'el\'eport 2 -- B.P. 179 \\
86960 Futuroscope CEDEX\quad France} \\
panyush@mathrs.sp2mi.univ-poitiers.fr}
}


\begin{thebibliography}{Pa95}
{\footnotesize
\bibitem[Al79]{andrei} {\sc Alekseevski\u i, A.V.}: The component groups of
the centralizers of unipotent elements in semisimple algebraic groups,
{\it Trudy Tbiliss. Matem. Inst. Akad. Nauk Gruzin. SSR} 
{\bf 62}(1979), 5--27 (Russian).

\bibitem[BK79]{bk} {\sc Borho, W., Kraft, H.}: \"Uber Bahnen und deren 
Deformationen bei Aktionen reduktiver Gruppen, {\it Comment. Math. Helv.} 
{\bf 54}(1979), 61--104.

\bibitem[CM93]{CoMc}
{\sc Collingwood, D.H.; McGovern, W.M.}: "Nilpotent orbits in semisimple
  Lie algebras", New York: Van Nostrand Reinhold, 1993.

\bibitem[Dy52]{EBD} 
{\sc Dynkin, E.B.}: Semisimple subalgebras of semisimple Lie algebras, 
{\it Matem. Sbornik} {\bf 30}(1952), {\rus N0}\,2, 349--462
(Russian). English translation: {\it Amer. Math. Soc. Transl.} II~Ser.,
{\bf 6}~(1957), 245--378.

\bibitem[El75]{El75}
{\sc Elashvili, A.G.}: The centralizers of nilpotent elements in semisimple Lie
algebras, {\it Trudy Tbiliss. Matem. Inst. Akad. Nauk Gruzin. SSR} 
{\bf 46}(1975), 109--132 (Russian).

\bibitem[El85]{El85}
{\sc Elashvili, A.G.}: Sheets in the exceptional Lie algebras, 
in: {\it ``Issledovaniya po Algebre''}, Tbilisi
1985, pp. 171--194 (Russian).

\bibitem[EP99]{wir}
{\sc Elashvili, A.G.; Panyushev, D.}:
Towards a classification of principal nilpotent pairs,
{\it Appendix} to [Gi99].

\bibitem[ElPa]{wir2}
{\sc Elashvili, A.G.; Panyushev, D.}:
A classification of principal and almost principal nilpotent pairs,
{\it in preparation}.

\bibitem[Gi99]{vitya} {\sc Ginzburg, V.}: Principal nilpotent pairs
in a semisimple Lie algebra I, {\it Preprint} math/9903059.

\bibitem[Ka82]{pasha} {\sc Katsylo, P.I.}: Sections of sheets in a reductive
algebraic Lie algebra, {\it Izv. AN SSSR. Ser. Matem.} {\bf 46}(1982),
477--486 (Russian). English translation: {\it Math. USSR-Izv.} 
{\bf 20}(1983), 449--458.

\bibitem[Ke83]{ke83} {\sc Kempken, G.}: 
Induced conjugacy classes in classical Lie-algebras,
{\it Abh. Math. Sem. Univ. Hamburg} {\bf 53}(1983), 53--83.

\bibitem[Ko59]{ko59} {\sc Kostant, B.}: The principal three-dimensional
subgroup and the Betti numbers of a complex simple Lie group, 
{\it Amer. J. Math.} {\bf 81}(1959), 973--1032.

\bibitem[Ko63]{bert} {\sc Kostant, B.}: Lie group representations in
polynomial rings, {\it Amer. J. Math.} {\bf 85}(1963), 327--404.

\bibitem[Pa99]{sph} {\sc Panyushev, D.}: On spherical nilpotent orbits
and beyond, {\it Ann. Inst. Fourier} {\bf 49}(1999), (to appear).

\bibitem[Ri79]{roger} {\sc Richardson, R.W.}: Commuting varieties of 
semisimple Lie algebras and algebraic groups, {\it Compositio math.} 
{\bf 38}(1979), 311--327.

\bibitem[Ru84]{nappes} {\sc Rubenthaler, H.}:  Param\'etrisation d'orbites
dans les nappes de Dixmier admissibles, {\it M\'emoire Soc. Math. France}, 
{\bf 15}(1984), 255--275.

\bibitem[Ru94]{rub} {\sc Rubenthaler, H.}: {\it ``Les paires duales dans les 
alg\`ebres de Lie reductives"}, Asterisque {\rus N0}\,219, S.M.F. 1994.

\bibitem[Sp82]{Sp} {\sc Spaltenstein. N.}: {\it ``Classes Unipotentes et 
Sous-groupes de Borel"}, Lecture notes in Math. {\bf 946},  Berlin 
Heidelberg New York: Springer 1982.

\bibitem[SpSt]{SS}
{\sc Springer, T.A.; Steinberg, R.}: Conjugacy classes, In: 
{\it "Seminar on algebraic groups and related finite groups".} Lecture notes 
in Math. {\bf 131}, pp.167--266,  Berlin 
Heidelberg New York: Springer 1970.

}
\end{thebibliography}
\end{document}